\documentclass[leqno,11pt]{article}
\usepackage{latexsym}
\usepackage{amssymb}
\usepackage{amsfonts}
\usepackage{amsmath}
\usepackage{epsfig, graphics, hyperref}

\makeatletter
\long\def\unmarkedfootnote#1{{\long\def\@makefntext##1{##1}\footnotetext{#1}}}
\makeatother

\newtheorem{definition}{Definition}[section]

\newtheorem{theorem}[definition]{Theorem}

\newtheorem{proposition}[definition]{Proposition}

\newtheorem{corollary}[definition]{Corollary}

\newtheorem{remark}[definition]{Remark}

\newtheorem{example}[definition]{Example}

\def\o{\Omega}

\def\mo{\hn (M )}

\def\m2{\hn (M ) /2}
\def\M2{\frac{\hn (\Omega )}{2}}
\def\u+{u_+^*}

\def\-p{\overline{p}}

\def\w0{{W_0^{m,A}(\Omega)}}

\def\R{\mathbb R}
\def\N{\mathbb N}
\def\Z{\mathbb Z}
\def\rn{{{\R}^n}}

\newcommand{\huno}{{\mathcal H}^{1}}
\newcommand{\hdue}{{\mathcal H}^{2}}
\newcommand{\hh}{{\mathcal H}^{n-1}}
\newcommand{\hn}{{\mathcal H}^{n}}
\newcommand{\medint}{-\kern  -,395cm\int}
\newcommand{\medintinrigo}{-\kern  -,315cm\int}
\newcommand{\medelle}{-\kern  -,235cm L}
\newcommand{\medellenrigo}{-\kern  -,180cm L}
\newcommand{\qed}{\thinspace\null\nobreak\hfill
\hbox{\vbox{\kern-.2pt\hrule height.2pt
depth.2pt\kern-.2pt\kern-.2pt \hbox to1.8mm {\kern-.2pt\vrule
width.4pt \kern-.2pt\raise1.8mm\vbox to.2pt{} \lower0pt\vtop
to.2pt{}\hfil\kern-.2pt \vrule
width.4pt\kern-.2pt}\kern-.2pt\kern-.2pt \hrule height.2pt
depth.2pt \kern-.2pt}}\par\medbreak}

\title{Bounds for eigenfunctions of the Laplacian \\ on noncompact Riemannian manifolds} 
\numberwithin{equation}{section}
\author{
  Andrea Cianchi\\
 {\it Dipartimento di Matematica e Applicazioni per
l'Architettura, Universit\`a di Firenze}\\ {\it Piazza Ghiberti
27, 50122 Firenze, Italy
}
\bigskip
 \and
 Vladimir G. Maz'ya \\
 {\it  Department of Mathematical Sciences, M\&O Building}\\ {\it University of Liverpool, Liverpool L69 3BX,
 UK}\\ and \\
{\it   Department of Mathematics, Link\"oping University, SE-581
83 Link\"oping, Sweden}}
\date{}

\begin{document}
\maketitle

\begin{abstract}\noindent
We deal with eigenvalue problems for the Laplacian on  noncompact
Riemannian manifolds $M$ of finite volume. Sharp conditions
 ensuring
 $L^q(M)$
 and $L^\infty (M)$ bounds for  eigenfunctions are
exhibited in terms of either the isoperimetric function  or the
isocapacitary function of $M$.
%

%

\end{abstract}


\unmarkedfootnote {
}

\section{Introduction }\label{sec1}

We are concerned with a class of eigenvalue problems for the
Laplacian on  $n$-dimensional Riemannian manifolds $M$ whose weak
formulation is:
\begin{equation}\label{weaksol}
\int_M \langle \nabla u \, , \nabla v \rangle \, d\hn (x)=\gamma
\int_M u\,v\, d\hn (x)
\end{equation}
for every test function $v$ in the Sobolev space $W^{1, 2}(M)$.
Here,  $u \in W^{1, 2}(M)$ is an eigenfunction associated with the
eigenvalue $\gamma \in \R$, $\nabla $ is the gradient operator,
 $\hn$ denotes the $n$-dimensional  Hausdorff measure on $M$,
i.e. the volume measure on $M$  induced by its Riemannian metric,
and $\langle \cdot \, , \cdot \rangle$ stands for the associated
scalar product.
\par\noindent Note that  various
special instances are included in this framework. For example, if
$M$ is a complete Riemannian manifold, then \eqref{weaksol} is
equivalent to a weak form of the equation
\begin{equation}\label{eigenproblem}
\Delta u + \gamma u =0 \quad \quad  \hbox{on $M$.}
\end{equation}
 In the case when $M$ is  an open subset of a
Riemannian manifold, and in particular of the Euclidean space
$\rn$, equation \eqref{weaksol} is a weak form of the eigenvalue
problem obtained on coupling equation \eqref{eigenproblem}
with homogeneous Neumann
boundary conditions.
\par
%
%
It is well known that quantitative information on eigenvalues and
eigenfunctions for elliptic  operators  in open subsets of
Euclidean space $\rn$ can be derived in terms of geometric
quantities associated with the domain.
 The quantitative analysis of spectral problems, especially for the
Laplacian,  on Riemannian manifolds is also very classical.
A great deal of contributions to
this topic regard  compact manifolds.
We do not even attempt to provide an exhaustive bibliography on
contributions to this matter; let us just mention the reference
monographs \cite{Chavel, BGM}, and the  papers \cite{Bou, Br, BD,
Cheeger, CGY, DS, Do1, Do2, Escobar,  Gallot, Gr, HSS, JMS, Na,
SS, So, SZ, Ya}.
%
\par
The present paper focuses the case when
$$
\hbox{$M$  need not be compact},$$ although
\begin{equation*}
\hn (M ) < \infty \,,
\end{equation*}
an assumption which will be kept in force throughout. We shall
also assume that $M$ is connected.
\par
 We are  concerned with estimates for  Lebesgue norms of eigenfunctions of the Laplacian on $M$.
When $M$ is compact, one easily infers, via local regularity
results for elliptic equations, that any eigenfunction $u$ of the
Laplacian belongs to $L^\infty (M)$. Explicit  bounds, with sharp
dependence on the eigenvalue $\gamma$, are also available
\cite{SS, SZ}, and require  sophisticated tools from differential
geometry and harmonic analysis.
If the compactness assumption is dropped, then the membership of
$u$ in $W^{1,2}(M)$ only (trivially) entails that $u \in L^2(M)$.
Higher integrability of eigenfunctions is not guaranteed anymore.
\par
Our aim is to exhibit minimal assumptions on $M$ ensuring $L^q(M)$
bounds for all $q<\infty$, or even $L^\infty (M)$ bounds for
eigenfunctions of the Laplacian on $M$. The results  that will be
presented can easily be extended to   linear uniformly elliptic
differential operators, in divergence form, with merely measurable
coefficients on $M$. However, we emphasize that our estimates are
new even for the Neumann Laplacian on open subsets of $\rn$ of
finite volume.
\par
The geometry of the manifold $M$ will come into play through
either the isocapacitary function $\nu _M$,  or the isoperimetric
function $\lambda _M$ of $M$. They are the largest functions of
the measure of subsets of $M$ which can be estimated by the
capacity, or by the perimeter of the relevant subsets.
%
%
%
Loosely speaking,
the asymptotic behavior of $\nu _M$ and $\lambda _M$ at $0$
accounts for the
 regularity of the geometry
 of the noncompact manifold $M$: decreasing this regularity  causes  $\nu _M (s)$ and $\lambda
_M (s)$ to decay faster to $0$ as $s$ goes to $0$.
\par
The  inequalities  associated with $\nu _M$ and $\lambda _M$ are
called the isocapacitary inequality and the isoperimetric
inequality on $M$, respectively. Thus, the isoperimetric
inequality on $M$ reads
\begin{equation}\label{isopineq}
\lambda _M( \hn (E)) \leq P(E)
\end{equation}
for every measurable set $E \subset M$ with $ \hn (E) \leq \m2$,
where
$\lambda _M: [0 , \m2 ] \rightarrow [0, \infty )$.
  \par\noindent In the isocapacitary
inequality that we are going to exploit,  the perimeter  on the
right-hand side of \eqref{isopineq} is replaced by the condenser
capacity of $E$ with respect to any subset $G \supset E$. The
resulting inequality has the form
\begin{equation}\label{isocineq}
\nu _M (\hn (E)) \leq C (E,G)
\end{equation}
for every measurable sets $E\subset G \subset M$ with $\hn (G)
\leq \m2$. Here,  $C(E,G)$ denotes the capacity of the condenser
$(E;G)$,
 and ${\nu _M}:
[0 , \m2 ] \rightarrow [0, \infty ]$
(see  Section \ref{sec2} for  precise definitions).
\par\noindent Introduced
 in \cite{Ma0}, the isoperimetric function $\lambda _M$ has been employed to provide necessary and sufficient conditions for embeddings
 of the Sobolev space $W^{1,1}(M)$ when $M$ is a domain in $\rn$
 \cite{Ma0}, and in a priori estimates for solutions to elliptic
 boundary value problems \cite{Ma0bis, Ma3}. Isocapacitary functions were introduced
and used in \cite{Ma0, Ma4, Ma5, Ma6, Ma2} in the characterization
of Sobolev embeddings for $W^{1,p}(M)$, with $p>1$, when $M$ is a
domain in $\rn$. Extensions to the case of Riemannian manifolds
can be found in \cite{Gr1, Gr}.
\par Both the conditions in
terms of  $\nu _M$, and those in terms of $\lambda _M$,  for
eigenfunction estimates in $L^q(M)$ or $L^\infty (M)$ that will be presented are sharp in
the class of  manifolds $M$ with prescribed asymptotic behavior of
$\nu _M$ and $\lambda _M$ at $0$. Each one of these two
approaches
has its own advantages. The isoperimetric function $\lambda _M$
has a transparent geometric character, and it is usually easier to
investigate. The isocapacitary function  can be less simple to
compute; however its use is in a sense more appropriate in the
present framework since it not only implies the results involving
$\lambda _M$, but leads  to finer conclusions in general.
Typically, this is the case when  manifolds with complicated
geometric configurations are taken into account.
\par
As for the proofs, let us just mention here that
crucial use is
made
of the isocapacitary inequality \eqref{isocineq} applied when $E$
is any level set of an eigenfunction $u$.
Note that
 customary methods, such as Moser iteration technique, which can be
exploited to derive eigenfunction estimates in classical
situations (see e.g. \cite{Saloff}), are  of no utility in the
present framework. In fact, Moser technique
would require a Sobolev embedding theorem for $W^{1,2}(M)$ into
some Lebesgue space smaller than $L^2(M)$, and this will  not be
guaranteed under the assumptions of our results.
%
%
%
%
%
%
%
\par
The  paper is organized as follows. The main results are stated in
the next section. The subsequent Section \ref{sec2} contains  some
basic definitions and properties concerning perimeter and capacity
which enter in our discussion.
In Section \ref{rev} we analyze a  class of manifolds of
revolution,
 which are used as model manifolds in the proof of the optimality of
 our results and in some examples. In particular, the behavior of their  isoperimetric
 and isocapacitary functions is investigated. Proofs of our bounds in $L^q(M)$  and
 in $L^\infty (M)$  are the object of Section
 \ref{sec5} and Section \ref{sec4}, respectively, where
explicit estimates depending on eigenvalues are also provided. The
final  Section \ref{appl} deals with  applications  of our
 results to two special instances: a family of manifolds of revolution
 whose profile has a borderline exponential behavior, and a family
 of  manifolds with  a sequence of clustering
 mushroom-shaped  submanifolds. In particular, the latter example  demonstrates that the use of $\nu _M$ instead of $\lambda _M$ can
actually lead to stronger results
when the regularity of eigenfunctions of the Laplacian is in
question.

\section{Main results }\label{main}
Our results will involve the manifold $M$ only through  the
asymptotic behavior of either $\nu_M$, or  $\lambda _M$ at $0$. They are stated in Subsections
\ref{maincap} and \ref{mainper}, respectively. \par Although the
criteria involving $\lambda _M$ admit independent proofs, along
the same lines as those involving $\nu _M$,
 the former will be deduced from the latter via
 the inequality:
\begin{equation}\label{lamnu}
\frac{1}{\nu _M(s)} \leq \int _{s}^{\m2 } \frac{dr}{\lambda _M
(r)^2} \qquad \hbox{for $s\in (0, \m2 )$,}
\end{equation}
which holds for any manifold $M$ (see the proof of
\cite[Proposition 4.3.4/1]{Ma2}).
 Let us notice that  a reverse inequality in \eqref{lamnu} does not hold in
general, even up to a multiplicative constant.
\par The proofs of the sharpness of the criteria for $\lambda _M$ and  $\nu _M$ require essentially the same  construction. We shall again
focus on the latter, and we shall explain how the relevant
construction also applies to the former.

\subsection{Eigenvalue estimates via  the isocapacitary
function of $M$}\label{maincap}


We begin with an optimal condition on the decay of $\nu _M$ at $0$
ensuring $L^q(M)$ estimates for eigenfunctions of the Laplacian on
$M$ for $q \in (2, \infty)$. Interestingly enough, such a condition is independent of $q$.

\begin{theorem}\label{eigencor} {\bf [$L^q$ bounds for eigenfunctions via $\nu _M$]}
Assume that
\begin{equation}\label{1001}
\lim _{s\to 0} \frac{s}{\nu _M(s)} =0\,.
\end{equation}
 Then for any
$q\in (2, \infty )$ and for any eigenvalue $\gamma$, there exists a constant $C= C(\nu _M , q,
\gamma )$ such that
\begin{equation}\label{ei5}
\|u \|_{L^q (M )} \leq C \|u \|_{L^{2} (M )}\,
\end{equation}
for every eigenfunction $u$ of the Laplacian on $M$ associated with $\gamma$.
\end{theorem}

An estimate for the constant $C$ in inequality \eqref{ei5}  can
also be provided -- see Proposition \ref{eigencorconstant},
Section \ref{sec5}.
\par
Let us note that condition \eqref{1001} turns out to be equivalent
to the compactness of the embedding $W^{1,2}(M) \to L^2(M)$
\cite{CM}. Hence, in particular, the variational characterization
of the eigenvalues of the Laplacian on $M$ entails that they
certainly exist under \eqref{1001}.
\par

Incidentally, let us also mention that, when $M$ is a complete
manifold, condition \eqref{1001} is also equivalent to the
discreteness of the spectrum of the Laplacian on $M$ \cite{CM}.

\smallskip
\par
The next result shows that assumption \eqref{1001} is essentially
minimal in Theorem \ref{eigencor}, in the sense that $L^q(M)$
regularity of eigenfunctions may fail under the mere assumption
that
$$ \nu _M (s) \approx s \qquad \quad \hbox{near $0$.}
$$
%
%
%
Here, and in what follows, the notation
\begin{equation}\label{approx}
f \approx g \qquad
 \end{equation}
 for functions $f,g : (0, \infty )
\to [0, \infty )$ means that there exist positive constants  $c_1$ and $c_2$ such
that
 \begin{equation}\label{approx1} c_1g(c_1s) \leq f(s) \leq c_2g(c_2s) \qquad \hbox{for $s>0$.}
 \end{equation}
Condition \eqref{approx1} is said to hold near $0$, or  near
infinity, if there exists a constant $s_0\,
>0$ such that \eqref{approx1} holds for $0 < s \leq
 s_0$ or for $s \geq s_0$, respectively.
 \par\noindent As in \eqref{1001}, all criteria  that will be presented are invariant under
 replacement of  $\nu _M$ or $\lambda _M$ with  functions $\approx$ near $0$.
 \begin{theorem}\label{Lqsharp}{\bf [Sharpness of condition \eqref{1001}]}
 For any $ n \geq 2$ and $q \in (2, \infty]$, there exists an $n$-dimensional Riemannian manifold $M$
such that
\begin{equation}\label{borderline}\nu _M (s) \approx s \qquad \quad \hbox{ near $0$},
\end{equation}
 and the Laplacian on $M$ has an eigenfunction $u \notin
L^q(M)$.
\end{theorem}

The important case when $q=\infty$, corresponding to the problem
of the boundedness of eigenfunctions, is not included in Theorem
\ref{eigencor}. This is the object of the following result, where
a slight strengthening of assumption \eqref{1001} is shown to
yield $L^\infty (M)$ estimates for eigenfunctions of the Laplacian
on $M$.

%
%

\begin{theorem}\label{eigencorinf} {\bf [Boundedness of eigenfunctions via $\nu _M$]}
Assume that
\begin{equation}\label{eiinftylambda}
\int _0 \frac {ds}{\nu _M (s)}  < \infty\,.
\end{equation}
Then for any eigenvalue $\gamma $, there exists a constant $C=C(\nu _M , \gamma )$
such that
\begin{equation}\label{ei5'}
\|u \|_{L^\infty (M )} \leq C
 \|u \|_{L^{2} (M )}\,
\end{equation}
for every eigenfunction $u$ of the Laplacian on $M$ associated with $\gamma$.
\end{theorem}

 An estimate for the constant $C$ in inequality \eqref{ei5'} is given in
 Proposition \ref{eigencorinfconstant}, Section \ref{sec4}.
\smallskip

Condition \eqref{eiinftylambda}   in Theorem \ref{eigencorinf} is
essentially sharp for the boundedness of eigenfunctions of the
Laplacian on $M$. In particular, it  cannot be relaxed to
\eqref{1001}, although the latter ensures $L^q(M)$ estimates for
every $q <\infty$.
 Indeed, under some mild
qualification,
Theorem \ref{eigensharp} below asserts  that given (up to
equivalence) any  isocapacitary function fulfilling \eqref{1001}
but not \eqref{eiinftylambda}, there exists a manifold $M$ with
the prescribed isocapacitary function on which the Laplacian has
an unbounded eigenfunction.
%
\par\noindent
A precise statement of this result involves the notion of function
of class $\Delta _2$. Recall that a non-decreasing function $f :
(0, \infty ) \to [0, \infty )$ is said to belong to the class
$\Delta _2$ near $0$ if there exist constants $c$ and $s_0$ such
that
\begin{equation}\label{delta2}
f(2s) \leq c f(s) \qquad \hbox{if $0 < s \leq s_0$.}
\end{equation}

\begin{theorem}\label{eigensharp}{\bf [Sharpness of condition \eqref{eiinftylambda}
]} Let $\nu$ be a non-decreasing function,
vanishing only at $0$, such that
\begin{equation*}
\lim _{s\to 0} \frac{s}{\nu (s)} =0\,,
\end{equation*}
but
\begin{equation*}
\int _0 \frac {ds}{\nu (s)}  =\infty\,.
\end{equation*}
Assume in addition that $\nu \in \Delta _2$ near $0$, and that
either $n \geq 3$ and
\begin{equation}\label{nu1tre}
\frac {\nu (s)}{s^{\frac{n-2}{n}}} \,\,\, \hbox{$\approx$
a
non-decreasing function near $0$,}
\end{equation}
or  $n=2$ and there exists $\alpha >0$ such that
\begin{equation}\label{nu1due}
\frac {\nu (s)}{s^{\alpha}} \,\,\, \hbox{
$\approx$
 a
non-decreasing function near $0$.}
\end{equation}
Then, there exists an $n$-dimensional Riemannian manifold $M$
fulfilling
\begin{equation}\label{1004}
\nu _M (s) \approx \nu (s)\qquad \hbox{near $0$,}
\end{equation}
 and such that
 the Laplacian on $M$ has
 an unbounded eigenfunction.
\end{theorem}
Assumption \eqref{nu1tre} or \eqref{nu1due} in Theorem
\ref{eigensharp} is explained by
the fact that, if
$M$ is compact, then
\begin{equation}\label{nucompact}
\nu _M(s) \approx
\begin{cases}
s^{\frac{n-2}{n}} &  \hbox{ if $n \geq 3$,}
\\
  \big( \log \frac{1}{s}\big)^{-1}
& \hbox{ if $n =2$,}
\end{cases}
 \end{equation}
near $0$, and that $\nu _M (s)$ cannot decay more slowly to $0$ as
$s\to 0$ in general. The assumption that $\nu \in \Delta _2$ near
$0$ is due to technical reasons.


\subsection{Eigenvalue estimates via the isoperimetric function
}\label{mainper}

The following criterion for $L^q(M)$ bounds of eigenfunctions
 in terms of the isoperimetric function $\lambda _M$ can be derived via  Theorem \ref{eigencor} and  inequality \eqref{lamnu}.
\smallskip
 \par
\par
\begin{theorem}\label{Lqlambdabis}{\bf [$L^q$ bounds for eigenfunctions via $\lambda _M$]}
Assume that
\begin{equation}\label{1001bis}
\lim _{s\to 0} \frac{s}{\lambda _M(s)} =0\,.
\end{equation}
Then for any $q\in (2, \infty )$ and any eigenvalue $\gamma$,
there exists a constant $C=C(\lambda _M , q, \gamma )$ such that
\begin{equation}\label{ei5lambda}
\|u \|_{L^q (M )} \leq C \|u \|_{L^{2} (M )}\,
\end{equation}
for every eigenfunction $u$ of the Laplacian on $M$ associated with $\gamma$.
\end{theorem}
%
%
%
%
An analogue of Theorem \ref{Lqsharp} on the minimality of
assumption \eqref{1001bis} in Theorem \ref{Lqlambdabis} is
contained in the the next result, showing that,  for every $q>2$,
eigenfunctions which do not belong to $L^q(M)$ may actually exist
when
$$ \lambda _M (s) \approx s \qquad \quad \hbox{near $0$.}
$$
\begin{theorem}\label{Lqsharplambda}{\bf [Sharpness of condition \eqref{1001bis}]}
 For any $ n \geq 2$ and $q \in (2, \infty]$, there exists an $n$-dimensional Riemannian manifold $M$
such that
\begin{equation}\label{borderlinelambda}
\lambda _M (s) \approx s \qquad \quad \hbox{ near $0$},
\end{equation}
 and the Laplacian on $M$ has an eigenfunction $u \notin
L^q(M)$.
\end{theorem}

A condition on $\lambda _M$, parallel to \eqref{eiinftylambda},
 ensuring the boundedness of
eigenfunctions of the Laplacian on $M$ follows from Theorem
\ref{eigencorinf} and inequality \eqref{lamnu}.

%

\begin{theorem}\label{boundlambda}{\bf [Boundedness of eigenfunctions via $\lambda_M$]}
Assume that
\begin{equation}\label{eiinftylambdabis}
\int _0 \frac {s}{\lambda _M (s)^2}\,ds  < \infty\,.
\end{equation}
Then for any eigenvalue $\gamma$, there exists a constant $C= C(\lambda _M , \gamma )$ such
that
\begin{equation}\label{ei5'lambda}
\|u \|_{L^\infty (M )} \leq C
 \|u \|_{L^{2} (M )}\,
\end{equation}
for every eigenfunction $u$ of the Laplacian on $M$ associated with $\gamma$.
\end{theorem}

Our last result tell us that the gap between condition
\eqref{eiinftylambdabis}, ensuring $L^\infty (M)$ bounds for
eigenfunctions, and condition \eqref{1001bis},  yielding $L^q (M)$
bounds for any $q <\infty$, cannot be essentially filled.

\begin{theorem}\label{eigensharplambda}{\bf [Sharpness of condition \eqref{eiinftylambdabis}
]} Let $\lambda $ be a non-decreasing function,
vanishing only at $0$, such that
\begin{equation*}
\lim _{s\to 0} \frac{s}{\lambda (s)} =0\,,
\end{equation*}
but
\begin{equation*}
\int _0 \frac {s}{\lambda (s)^2} \,ds =\infty\,.
\end{equation*}
Assume in addition that
\begin{equation}\label{nu1trebis}
\frac {\lambda (s)}{s^{\frac{n-1}{n}}} \,\,\, \hbox{
$\approx$
a non-decreasing function near $0$.}
\end{equation}
Then, there exists an $n$-dimensional Riemannian manifold $M$
fulfilling
\begin{equation}\label{1004bis}
\lambda _M (s) \approx \lambda (s)\qquad \hbox{near $0$,}
\end{equation}
 and such that
 the Laplacian on $M$ has
 an unbounded eigenfunction.
\end{theorem}

Assumption \eqref{nu1trebis}  in Theorem \ref{eigensharplambda} is
required in the light of the fact that
\begin{equation}\label{lambdacompact}
\lambda _M(s) \approx s^{\frac{n-1}n} \qquad \quad \hbox{near $0$}
\end{equation}
for any compact manifold $M$, and that $\lambda _M (s)$ cannot
decay more slowly to $0$ as $s \to 0$
 in the noncompact case.

\section{Background and preliminaries}\label{sec2}

Let $E$ be a measurable subset of $M$. The perimeter $P(E)$ of $E$
is defined as
$$P(E) = \hh (\partial ^* E)\,,$$
where $\partial ^* E$ stands for the essential boundary of $E$ in
the sense of geometric measure theory, and
 $\hh$ denotes the $(n-1)$-dimensional Hausdorff measure on $M$ induced by its Riemannian metric
(\cite{AFP, Ma2}). Recall that $\partial ^* E$ agrees with
 the topological boundary $\partial E$ of $E$ when $E$ is
 sufficiently regular, for instance an open subset  of $M$ with a
 smooth boundary. In the special case when $M = \Omega$, an open subset
 of $\rn$,
 and $E \subset \Omega$, we have that $P(E)=\hh (\partial ^*_{\rn} E \cap
 \Omega)$, where $\partial ^*_{\rn} E$ denotes the essential
 boundary of $E$ in $\rn$.
\par\noindent
The isoperimetric function $\lambda _M$ of $M$ is defined as
\begin{equation}\label{isopfunct}
\lambda _M(s) = \inf \{P(E): s \leq \hn (E) \leq \m2 \} \qquad
\quad \hbox{for $s \in[0, \m2 ]$\,.}
 \end{equation}
The  isoperimetric inequality \eqref{isopineq}
 is just a rephrasing of definition \eqref{isopfunct}. The point
 is thus  to derive information about the function $\lambda _M$,
 which is explicitly known only for
 Euclidean balls and spheres \cite{BuZa, Cpoincare, Ma2}, convex cones \cite{LP},
and manifolds in special classes \cite{BC, CF, CGL, GP, HHM, Kl,
MJ, Pi, Ri}. Various qualitative and quantitative properties of
$\lambda _M$ are however available -- see e.g.
 \cite{BuZa,
Crelative, HK, KM, La, Ma2}.
%
%
In particular, since we are assuming that $M$ is connected,
 \begin{equation}\label{lambdapos}
 \lambda _M(s) >0 \qquad  \qquad   \hbox{for $s \in (0, \m2]$,}
 \end{equation}
as an analogous argument as in \cite[Lemma 3.2.4]{Ma2} shows.
\par
The Sobolev space $W^{1,p}(M )$ is defined, for $ p \in [1, \infty
]$, as
\begin{equation*}
W^{1,p}(M ) = \{u \in L^p(M): \,\,\hbox{$u$ is weakly
differentiable on $M$ and $|\nabla u| \in L^p(M )$ }\}\,.
\end{equation*}
We denote by $W^{1,p}_0(M)$  the closure in $W^{1,p}(M)$ of the
set of smooth compactly supported functions on $M$.
\par
The standard $p$-capacity of a set $E\subset M$ can be defined  as
\begin{equation}\label{-506}
C_p(E) = \inf \left\{ \int _M |\nabla u|^p\,dx : u\in
W^{1,p}_0(M), u \geq 1 \,\,\, \hbox{in some neighbourhood
of}\,\,\, E\right\}.
\end{equation}
A property concerning the pointwise behavior of functions is said
to hold $C_p$-quasi everywhere in $M$, $C_p$-q.e. for short,  if
it is fulfilled outside a set of $p$-capacity zero.
\par\noindent
Each function $u \in W^{1,p}(M)$ has a representative
$\widetilde{u}$,  called the precise representative, which is
$C_p$-quasi continuous, in the sense that for every $\varepsilon
>0$, there exists a set $A \subset M$, with $C_p(A) <
\varepsilon$, such that $f_{|M \setminus A}$ is continuous in $M
\setminus A$. The function $\widetilde{u}$ is unique, up to
subsets of $p$-capacity zero. In what follows, we assume that any
function $u \in W^{1,p}(M)$ agrees with its precise
representative.
\par
In the light of a classical result in the theory of capacity 
(\cite[Proposition 12.4]{Da}, \cite[Corollary 2.25]{MZ}), we adopt
the following definition of capacity of a condenser. Given sets $E
\subset G \subset M$, the capacity $C_p(E,G)$  of the condenser
$(E, G)$ relative to $\o$ is defined as
\begin{equation}\label{-508}
C_p(E,G) = \inf \left\{ \int _M |\nabla u|^p\,dx : u\in
W^{1,p}(M), \hbox{$u \geq 1$ $C_p$-q.e. in $E$ and $u \leq 0$
$C_p$-q.e. in $M \setminus G$} \right\}\,.
\end{equation}
\par\noindent
Accordingly, the $p$-isocapacitary function $\nu _{M, p}: [0, \m2
] \to [0, \infty ]$ of $M $ is  given by
\begin{multline}\label{-509}
\nu _{M, p}(s) = \inf \left\{C_p(E,G): \hbox{$E$ and $G$ are
measurable subsets of $M$ such that}\right.\\
\left. \hbox{$E \subset G \subset M$ and  $s \leq \hn (E) \leq \hn
(G) \leq \m2$} \right\}\qquad \hbox{for $s\in [0, \m2 ]$}.
\end{multline}
The function $\nu _{M,p}$ is clearly non-decreasing.  In what
follows, we shall always deal with the left-continuous
representative of $\nu _{M, p}$, which, owing to the monotonicity
of $\nu _{M, p}$, is pointwise dominated by the right-hand side of
\eqref{-509}. Note that
%
%
 \begin{equation}\label{nu1lambda}
 \nu _{M,1}= \lambda _M
 \end{equation}
   as shown  by an analogous argument as in
\cite[Lemma 2.2.5]{Ma2}.
\par
When $p=2$, the case of main interest in the present paper, we
drop the index $p$  in $C_p(E,G)$ and $\nu _{M,p}$, and simply set
$$ C(E,G)=C_2(E,G)\,,$$
and
$$ \nu _M = \nu _{M,2}.$$
By \eqref{lambdapos} and \eqref{lamnu}, one has that
\begin{equation}\label{nupos}
 \nu _M(s)> 0 \qquad \qquad \hbox{for $s \in [0, \m2 ]$.}
 \end{equation}
\par
For any measurable function $u$ on $M$, we define its distribution
function
 $\mu _u: \R \to
[0, \infty )$ as
$$\mu  _u(t) = \hn (\{x\in M : u(x) \geq t\}) \qquad \quad \hbox{for
$t \in \R$}.$$ Note that here $\mu _u$ is defined in terms of $u$,
and not of $|u|$ as customary. The signed decreasing rearrangement
$u^\circ : [0, \mo ] \to [-\infty , \infty ]$ of $u$ is given by
$$u^\circ (s) = \sup \{t : \mu _u(t) \geq s\} \qquad \quad \hbox{for
$s\in [0, \mo ]$.}$$ The median of $u$ is defined by
\begin{align}\label{mediana}
{\rm med } (u) = u^\circ (\m2 )\,.
\end{align}
Since $u$ and $u^\circ$ are equimeasurable functions, one has that
\begin{equation}\label{Lebesgue}
\|u^\circ \|_{L^q(0, \mo )} = \|u\|_{L^q(M)}
\end{equation}
for every $q \in [1, \infty ]$. Moreover, by an analogous argument
as in \cite[Lemma 6.6]{CEG}, if $u \in W^{1,p}(M)$ for some $p \in
[1, \infty ]$, then
\begin{equation}\label{ac}
\hbox{ $u^\circ $ is locally absolutely continuous in
    $(0, \mo)$.}
    \end{equation}
\par
Given $u \in W^{1,2}(M)$, we define
  the function $\psi _u : \R \to [0, \infty )$  as
 \begin{equation}\label{psi}
 \psi _u(t) = \int _0^t \frac{d\tau}
 {\int _{\{u=\tau \}}
|\nabla u|\,d \hh (x)} \qquad \quad \hbox{for $t\in \R$}\,.
\end{equation}
On making use of (a version on manifolds) of \cite[Lemma
2.2.2/1]{Ma2}, one can easily show that if
\begin{equation}\label{e-1}
{\rm med } (u) = 0\,,
\end{equation}
 then
\begin{equation}\label{e0'}
\nu _{M}(\hn (\{u\geq t\})) \leq \frac 1{\psi_u (t)} \quad
\hbox{for $t> 0$,}
\end{equation}
and
\begin{equation}\label{e0}
\nu _{M}(s) \leq \frac 1 {\psi _{u}(u^\circ (s))} \quad \qquad
\qquad\hbox{for $s\in (0, \m2 )$}.
\end{equation}

\section{Manifolds of revolution}\label{rev}

In this section we focus on  a class of manifolds of revolution to
be employed in our proofs of Theorems \ref{Lqsharp} and
\ref{eigensharp}. Specifically, we investigate on their
isoperimetric and isocapacitary functions.
\par
Let $L \in (0, \infty ]$, and let $\varphi : [0, L ) \to [0,
\infty )$ be a function in $C^1([0, L))$,
 such that
 \begin{equation}\label{n0}
\hbox{$\varphi (r) >0$ \quad  for  $r\in (0, L)$,}
\end{equation}
\begin{equation}\label{n1'}
\varphi (0) = 0\,, \qquad\hbox{and} \qquad  \varphi '(0) =1\,.
%
\end{equation}
Here, $\varphi '$ denotes the derivative of $\varphi$.
For $n \geq 2$, we call $n$-dimensional manifold of revolution $M$
associated with $\varphi$ the ball in $\mathbb R ^n$ given,
 in polar coordinates, by $\{(r, \omega ): r\in [0, L) ,
 \omega \in \mathbb S ^{n-1}\}$ and
 endowed with the Riemannian metric
\begin{equation}\label{metric}
ds^2 = d r^2 + \varphi (r)^2 d\omega ^2\,.
\end{equation}
Here,  $d \omega ^2$ stands for the standard metric on $ \mathbb S
^{n-1}$.
Owing to our assumptions on $\varphi$, the metric \eqref{metric}
is of
class $C^1(M)$. 
Note that, in particular,
\begin{equation}\label{n4''}
\int _M u \,d\hn = \int _{\mathbb S ^{n-1}} \int _0^L u \,\varphi
(r)^{n-1}\, dr\, d \sigma _{n-1}\,,
\end{equation}
for any integrable function $u : M \to \R$. Here, $\sigma _{n-1}$
denotes the $(n-1)$-dimensional Hausdorff measure on $\mathbb S
^{n-1}$.
\par
The length of the gradient of a function $u: M \to \mathbb R$ is
defined
 by $|\nabla u | = \sqrt
{\langle \nabla u \, , \nabla u \rangle}$, and takes the form
\begin{equation}\label{n11}
|\nabla u| = \sqrt{\Big(\frac{\partial u}{\partial r}\Big)^2 +
\frac 1{\varphi (r)^2} |\nabla _{\mathbb S ^{n-1}} u| ^2},
\end{equation}
where $\nabla _{\mathbb S ^{n-1}}$ denotes the gradient operator
on  $\mathbb S ^{n-1}$. Moreover, if $u$ depends only on $r$, then
\begin{equation}\label{n12}
\Delta u = \frac 1{\varphi (r)^{n-1}}\frac{d}{dr}\bigg(\varphi
(r)^{n-1}\frac{du}{dr}\bigg)\,.
\end{equation}
Thus, for functions $u$ depending only on $r$, equation
\eqref{eigenproblem} reduces to the ordinary differential equation
%
%
\begin{equation}\label{n106}
\frac{d}{dr}\bigg(\varphi (r)^{n-1}\frac{du}{dr}\Big) + \gamma
\varphi (r)^{n-1} u=0 \qquad \hbox{for $r \in (0, L)$.}
\end{equation}
The membership of $u$ in $W^{1,2}(M)$ reads
\begin{equation}\label{n107}
\int _0^L \bigg(u^2 +\bigg(\frac{du}{dr}\bigg)^2\bigg)\varphi
(r)^{n-1} dr < \infty\,.
\end{equation}
Now, fix any $r_0 \in (0, L)$, set
\begin{equation}\label{*}
s_0 = \int _{r_0}^L \frac {d\rho}{\varphi (\rho)^{n-1}},
\end{equation}
%
and define $\psi : (0, L) \to \R$ as
\begin{equation}\label{n108}
\psi (r) = \int _{r_0}^r \frac {d\rho}{\varphi (\rho)^{n-1}}
\qquad \hbox{for $r \in (0, L)$.}
\end{equation}
Under the change of variables
$$s = \psi (r),$$
$$v(s) = u(\psi ^{-1}(s)),$$
and
$$p(s) = \varphi \big(\psi ^{-1}(s)\big)^{2(n-1)},$$
equations \eqref{n106} and \eqref{n107} turn into
\begin{equation}\label{n109}
\frac{d^2v}{ds^2} + \gamma p(s) v =0\, \qquad \quad \hbox{for $s
\in (-\infty , s_0)$,}
\end{equation}
and
\begin{equation}\label{n111}
\int _{-\infty}^{s_0} \bigg(v^2 p(s) + \bigg(\frac{dv}{ds}\bigg)^2
\bigg)ds < \infty\,,
\end{equation}
respectively. \par\noindent
We finally note that if  $r \in (0, L)$
and $B(r)=\{(\rho, \omega ): \rho\in [0, r) ,
 \omega \in \mathbb S ^{n-1}\}$, a ball on $M$ centered at $0$, then
\begin{equation}\label{n6}
\hh (\partial (M\setminus B(r))) = \hh (\partial B(r)) = \omega
_{n-1} \varphi (r)^{n-1}\,,
\end{equation}
and
\begin{equation}\label{n7}
\hn (M\setminus B(r)) = \omega _{n-1} \int _r^L \varphi
(\rho)^{n-1}d\rho\,,
\end{equation}
where $\omega _{n-1} = \hh (\mathbb S ^{n-1})$.
\par
 The main result of this section is contained in the
following theorem, which provides us (up to equivalence) with the
functions $\lambda _M$ and $\nu _M$  for a manifold of revolution
$M$ as above. In what follows, we set $n'=\frac n{n-1}$, the
H\"older conjugate of $n$.
\begin{theorem}\label{disisopA}
Let $L\in (0, \infty ]$ and let $\varphi : [0, L) \to [0, \infty
)$ be a function in $C^1([0, L))$ fulfilling \eqref{n0} and
\eqref{n1'} and such that:
\par \noindent
\par \noindent (i) $\lim _{r \to L}\varphi (r)=0$;
\par\noindent  (ii) there exists $L_0 \in
(0, L)$ such that $\varphi$ is decreasing and convex in $(L_0,
L)$;
\par\noindent (iii) $\int _0^L \varphi (\rho)^{n-1} \, d\rho < \infty$.
\par\noindent
Then the metric of the $n$-dimensional manifold of revolution $M$
built upon $\varphi$ is of class $C^1(M)$, and $\hn (M) < \infty$.
Moreover, let $\lambda _0$ be the function implicitly defined by
\begin{equation}\label{l7bis}
\lambda _0 \bigg(\omega _{n-1}\int _r^L \varphi (\rho )^{n-1}d\rho
\bigg)= \omega _{n-1}\varphi (r )^{n-1} \qquad \hbox{for $r\in
(L_0, L)$}\,,
\end{equation}
and  such that $\lambda _0 (s) = \lambda _0 \big(\omega _{n-1}\int
_{L_0}^L \varphi (\rho )^{n-1}d\rho\big)$ for $s \in \big(0,
\omega _{n-1} \int _{L_0}^L \varphi (r)^{n-1}dr\big)$. Then
\begin{equation}\label{n50}
\lambda _M (s)\approx \lambda _0 (s) \qquad \hbox{near $0$}\,,
\end{equation}
and
\begin{equation}\label{n51}
\nu _M (s) \approx \frac 1{\int _s^{\m2}\frac{dr}{\lambda _0
(r)^2}}\qquad \hbox{near  $0$.}
\end{equation}
\end{theorem}

\par\noindent
{\bf Proof
} The fact that $M$ is a Riemannian manifold of class $C^1$
follows from assumptions \eqref{n0} and \eqref{n1'}. Furthermore,
by  (iii),
$$\hn (M) = \omega _{n-1}\int _0^L \varphi (\rho
)^{n-1}d\rho < \infty .$$
\par\noindent
As for \eqref{n50} and \eqref{n51}, let us begin by observing
that, since $\varphi$ is decreasing in $(L_0, L)$, the function
$\lambda _0$ is increasing in $\big(0, \omega _{n-1} \int _{L_0}^L
\varphi (r)^{n-1}dr\big)$. Moreover, there exists a constant $C$
such that
\begin{equation}\label{l52}
\lambda _0 (s) \leq C s^{1/n'}
\end{equation}
for $s \in \big(0, \omega _{n-1} \int _{L_0}^L \varphi
(r)^{n-1}dr\big)$. Indeed, since $\lim _{r\to L}\varphi (r)=0$ and
$-\varphi '(r)$ is a nonnegative non-increasing function in $(L_0,
L)$, one has that
\begin{align}\label{l53}
-\varphi '(L_0) \int _r^L \varphi (\varrho )^{n-1} d\varrho & \geq
-\varphi '(r) \int _r^L\varphi (\varrho )^{n-1} d\varrho \\
\nonumber & \geq  \int _r^L -\varphi '(\varrho)\varphi (\varrho
)^{n-1} d\varrho = \frac 1n \,\varphi (r)^n \qquad \hbox{for $r
\in (L_0, L)$,}
\end{align}
 whence \eqref{l52} follows, owing to
\eqref{l7bis}.
\par\noindent
Define the map $\Phi : M \setminus \overline{B(L_0)} \to \rn$
as
\begin{equation}\label{l12}
\Phi (r,\omega ) = (\varphi (r), \omega ) \qquad \hbox{for $(r,
\omega ) \in ({L_0}, L) \times \mathbb S ^{n-1}$.}
\end{equation}
Clearly, $\Phi $ is a diffeomorphism between $M \setminus
\overline{B(L_0)}$ and $\Phi \big(M \setminus
\overline{B(L_0)}\big)$.
 \par\noindent Given any smooth function $v : M \setminus
\overline{B({L_0})} \to \R$,
we have that
\begin{align}\label{l13}
\int _{\Phi (M \setminus \overline{B({L_0})})}&|\nabla (v \circ
\Phi ^{-1})| dx = \int _{{\mathbb S} ^{n-1}} \int _{0}^{\varphi
({L_0})} \sqrt{(v \circ \Phi ^{-1})_\varrho ^2+ \frac 1{\varrho
^2}|\nabla _{\mathbb S ^{n-1}} (v \circ \Phi
^{-1})|^2}\,\varrho ^{n-1} d\varrho d\sigma _{n-1}\\
\nonumber & =
 \int _{\mathbb
S ^{n-1}} \int _{{L_0}}^L \sqrt{\frac{1}{\varphi '(r)^2}
\Big(\frac{\partial v}{\partial r}\Big)^2 + \frac 1{\varphi (r)
^2}|\nabla _{\mathbb S ^{n-1}} v
|^2}\,\varphi (r) ^{n-1}|\varphi '(r) | dr d\sigma _{n-1}
\\ \nonumber & \leq
2(1+ \sup _{r \in [L_0 , L)}|\varphi '(r)|) \int _{\mathbb S
^{n-1}} \int _{{L_0}}^L \sqrt{\Big(\frac{\partial v}{\partial
r}\Big)^2 + \frac 1{\varphi (r) ^2}|\nabla _{\mathbb S ^{n-1}} v
|^2}\,\varphi (r) ^{n-1} dr d\sigma _{n-1}
\\ \nonumber & = 2(1 - \varphi '(L_0))\int _{M\setminus \overline{B({L_0})}} |\nabla v| d\hn\,.
\end{align}
By approximation, the inequality between the leftmost side and the
rightmost side of \eqref{l13} continues to hold  for any function
of bounded variation $v$, provided that the integrals of the
gradients are replaced by the total variations. In particular, on
applying the resulting inequality to
 characteristic function of sets, we
obtain that
\begin{equation}\label{l14}
\hh (\partial (\Phi (E))) \leq C \hh (\partial E)\,
\end{equation}
for every smooth set $E\subset  M \setminus \overline{B({L_0})}$,
where $C=2(1 - \varphi '(L_0))$. Given any such set $E$, the
classical isoperimetric inequality in $\rn$ tells us that
\begin{equation}\label{l15}
n ^{1/n'}\omega _{n-1}^{1/n} \mathcal L ^n(\Phi (E))^{1/n'} \leq
\hh (\partial (\Phi (E)))\,,
\end{equation}
where $\mathcal L ^n$ denotes the Lebesgue measure in $\rn$.
 On the other hand,
\begin{equation}\label{l16}
\mathcal L ^n(\Phi (E)) = \int _{\mathbb S ^{n-1}} \int
_{0}^{\varphi (L_0)}\chi_{\Phi (E)} \varrho ^{n-1} d\varrho
d\sigma _{n-1} = \int _{\mathbb S ^{n-1}} \int _{L_0}^L \chi_{E}
\varphi (r)^{n-1} |\varphi '(r)| dr d\sigma _{n-1}\,,
\end{equation}
where $\chi _E$ and $\chi_{\Phi (E)}$ stand for the characteristic
functions of the sets $E$ and $\Phi (E)$, respectively.
 \par\noindent
 Define $\Lambda : [0, L) \to [0, \hn (M)]$ as
$$\Lambda (r) = \omega _{n-1}\int _r^L \varphi (\rho
)^{n-1}d\rho \qquad \quad \hbox{for $r \in [0, L)$,}$$
 whence, by \eqref{n7}, $\Lambda (r)= \hn (M \setminus B(r)$ for $r \in [0,
 L)$.
%
%
%
Since $|\varphi '|=\,-\varphi '$ in $(L_0 , L)$, and $- \varphi '$
 is a non-increasing function in $(L_0 , L)$,
\begin{multline}\label{l17}
\int _{\mathbb S ^{n-1}} \int _{L_0}^L \chi_{E} \varphi (r)^{n-1}
|\varphi '(r)| dr d\sigma _{n-1} \geq \int _{\mathbb S ^{n-1}}
\int _{\Lambda ^{-1}(\hn(E))}^L  \varphi (r)^{n-1} |\varphi '(r)|
dr d\sigma _{n-1} \\ = \omega _{n-1} \int _{\Lambda
^{-1}(\hn(E))}^L \varphi (r)^{n-1} (-\varphi '(r)) dr =
\frac{\omega _{n-1}}n \varphi (\Lambda ^{-1}(\hn(E)))^n\,.
\end{multline}
Combining \eqref{l14}-\eqref{l17} yields
\begin{equation}\label{l18}
C \varphi (\Lambda ^{-1}(\hn(E)))^{n-1} \leq \hh(\partial E)
\end{equation}
for some positive constant $C$. By \eqref{l7bis},
\begin{equation}\label{l19}
\omega _n \varphi (\Lambda ^{-1}(\hn(E)))^{n-1} =\lambda _0
\bigg(\omega _{n-1}\int _{\Lambda ^{-1}(\hn(E))}^L \varphi (\rho
)^{n-1}d\rho \bigg) = \lambda _0 (\hn(E))\,.
\end{equation}
From \eqref{l18} and \eqref{l19} we obtain that
\begin{equation}\label{l20}
C \lambda _0 (\hn(E)) \leq \hh(\partial E)\,,
\end{equation}
from some constant $C$ and any smooth set $E\subset M \setminus
\overline{B({L_0})}$.
%
\par
Now, let $L_1 \in (L_0, L)$ be such that $\hn (B(L_1)) > \m2$.
Observe that $\overline{B({L_1})}$ is a smooth
compact Riemannian submanifold of $M$  with boundary $\partial
B({L_1})$ diffeormorphic to a closed ball in $\rn$. Thus, an isoperimetric inequality of the form
\begin{equation}\label{l21}
C \hn(E)^{1/n'} \leq \hh(\partial E)
\end{equation}
holds for some constant $C$ and  for any set of finite perimeter
$E \subset \overline{B({L_1})}$. Moreover, there exists a positive
constant $C$  such that
\begin{equation}\label{l22}
\hh(E \cap \partial B( {L_1})) \leq C \hh(\partial E \cap B(
{L_1}))
\end{equation}
for any smooth set $E \subset \overline{B({L_1})}$ such that $\hn
(E) \leq \m2$ ($< \hn \big(\overline {B({L_1})}\big)$, by our
choice of $L_1$).
\par
Owing to \eqref{l20}-\eqref{l22}, for  any smooth set $E \subset
M$ such that $\hn  (E) \leq \hn (M)/{2}$
\begin{align}\label{l23}
\hh(\partial E) & = \hh(\partial (E \cap \overline{B( {L_1})})) +
\hh(\partial (E\cap (M\setminus B( {L_1})))) - 2\hh(E \cap
\partial B( {L_1}))\\ \nonumber
& \geq C \hn(E \cap \overline{B( {L_1})})^{1/n'} + C \lambda _0
(\hn(E\cap (M\setminus B( {L_1}))) - C \hh( B( {L_1}) \cap
\partial E)
\\ \nonumber
& \geq C \hn(E \cap \overline{B( {L_1})})^{1/n'} + C \lambda _0
(\hn(E\cap (M\setminus B( {L_1}))) - C \hh(
\partial E)\,,
\end{align}
for some positive constant $C$.  Consequently, there exists a
constant $C$ such that
\begin{align}\label{l24}
C \hh(\partial E) \geq \hn(E \cap \overline{B( {L_1})})^{1/n'} +
 \lambda _0 (\hn (E\cap (M\setminus B( {L_1})))
 \end{align}
 for  any smooth set $E \subset
M$ such that $\hn  (E) \leq \hn (M)/{2}$. Now, we claim that there
exists a constant $C$ such that
 such that
 \begin{equation}\label{l25}
 s^{1/n'} + \lambda _0 (\sigma ) \geq C \lambda _0 \Big(\frac{s+\sigma
 }{2}\Big) \qquad \quad \hbox{for $s, \sigma \in [0, \hn (M)/{2}]$.}
 \end{equation}
Indeed, if $\sigma \leq s$, then, by \eqref{l52},
\begin{equation}\label{l26}
 s^{1/n'}  \geq C \lambda _0 (s) \geq C \lambda _0 \Big(\frac{s+\sigma
 }{2}\Big)
 \end{equation}
 for some positive constant $C$, whereas, if $s \leq \sigma$, then trivially
\begin{equation}\label{l27}
 \lambda _0 (\sigma )  \geq  \lambda _0 \Big(\frac{s+\sigma
 }{2}\Big)\,.
 \end{equation}
 Coupling \eqref{l24} and \eqref{l25} yields
\begin{align}\label{l28}
C \lambda _0 (\hn (E)/2) \leq  \hh(\partial E)
 \end{align}
 for some positive constant $C$ and for every smooth set $E \subset M$ such that $\hn  (E) \leq \hn (M) /{2}$.
 By approximation, inequality \eqref{l28} continues to hold for every set of
 finite perimeter $E \subset M$ such that $\hn  (E) \leq \hn (M) /{2}$.
  Since $\lambda _0 $ is non-decreasing, inequality \eqref{l28} ensures
 that
\begin{equation}\label{l29}
{\lambda} _M (s) \geq C \lambda _0 (s/2) \qquad \hbox{for $s\in
[0, \m2]$.}
\end{equation}
 \par\noindent
 On the other hand, equation \eqref{l7bis} entails that ${\lambda} _M (s)\leq \lambda _0 (s)
 $ for small $s$, and hence there exists a constant $C$ such that
 \begin{equation}\label{l30}
 {\lambda} _M (s)\leq C \lambda _0 (s) \qquad \hbox{for $s \in [0 , \hn (M)/2]$.}
 \end{equation}
Equation \eqref{n50} is fully proved.
\par
As far as \eqref{n51} is concerned, by \eqref{lamnu} and
\eqref{l29},
 \begin{align}\label{nuotto}
\frac 1{\nu _M(s)} & \leq \int _s^{\hn (M)/2} \frac {dr}{{\lambda}
_M(r)^2} \leq \frac 1{C^2}\int _s^{\hn (M)/2} \frac {dr}{\lambda
_0 (r/2)^2} 
\\
\nonumber & \leq \frac 2{C^2}\int _{s/2}^{\hn (M)/2} \frac
{dr}{\lambda _0 (r)^2} \qquad \hbox{for $s \in (0 , \hn (M)/2]$.}
\end{align}
\par\noindent
In order to prove a reverse inequality, set $R = \max\{L_0,
\Lambda ^{-1}(\hn (M)/2)\}$.
Moreover, given $s \in (0, \m2 )$, let $r \in (R, L)$ be such that
\begin{equation}\label{nu9}
s = \hn (M \setminus B(r)) = \omega _{n-1}\int _r ^L \varphi (\tau
)^{n-1} \,d\tau\,.
\end{equation}
Let $u=u(\rho )$ be the function given by
$$ u(\rho ) = \begin{cases}
0 & \hbox{if $\rho \in (0,R]$,} \\
\displaystyle{\frac{\int _{R} ^\rho \frac {d\tau}{\varphi (\tau
)^{n-1}}}{\int _R ^r \frac {d\tau}{\varphi (\tau )^{n-1}}} }&
\hbox{if $\rho \in (R,
r)$,}\\
1 & \hbox{if $\rho \in [r,L)$,}
\end{cases}
$$
and let
$$ E= M \setminus B(r)\qquad \hbox{and}\qquad G=M \setminus
B(R).$$
%
Hence,
\begin{align}\label{nu10}
\nu _M(s)  & \leq C(E,G) \leq \int _M |\nabla u |^2 d \hn (x) =
 \displaystyle{\frac {\omega
_{n-1}}{\int _{R} ^r \frac {d\tau}{\varphi (\tau )^{n-1}}}} \\
\nonumber & = \displaystyle{\frac 1{\int _s^{\Lambda (R)} \frac
{d\rho }{\lambda _0 (\rho )^2}}} \leq \displaystyle{\frac C{\int
_s^{\hn (M)/2} \frac {d\rho }{\lambda _0 (\rho )^2}}}
\qquad \qquad \quad \hbox{for $s \in (0 , \hn (M)/2]$,}
\end{align}
for  some constant $C$.  Note that the second equality is a
consequence of \eqref{l7bis}, owing to a change of variable.
Equation \eqref{n51}  follows from \eqref{nuotto} and
\eqref{nu10}.
 \qed
From  \eqref{n50} and \eqref{n51}, it is easily verified that
conditions \eqref{1001} and \eqref{1001bis} are equivalent for
manifolds of revolution as in Theorem \ref{disisopA}. Moreover,
 these conditions can be characterized in terms of the function $\varphi$ appearing in its statement. The
same observation applies to \eqref{eiinftylambda} and
\eqref{eiinftylambdabis}. These observations  are summarized in
the following statement.
\begin{corollary}\label{manifoldsofrevolution}
Let $\varphi$ be as in Theorem \ref{disisopA}, and let $M$ be the
$n$-dimensional manifold of revolution built upon $\varphi$. Then:
\par\noindent
(i) Conditions \eqref{1001}, \eqref{1001bis},   and
\begin{equation*}
\lim _{r\to L } \bigg(\int _R^r\frac {d\varrho}{\varphi
(\varrho)^{n-1}}\bigg)\bigg(\int _r^L \varphi (\varrho
)^{n-1}d\varrho\bigg) =0\,
\end{equation*}
%
%
for any $R \in (0, L)$ are equivalent.
\par\noindent
(ii) Conditions
 \eqref{eiinftylambda},   \eqref{eiinftylambdabis}, and
\begin{equation*}
\int ^L \bigg(\frac 1{\varphi (r)^{n-1}}\int _r^L \varphi (\rho
)^{n-1}d\rho\,\bigg)\, dr < \infty
\end{equation*}
are equivalent.
\end{corollary}

The remaining part of this section is devoted to showing that,
given functions $\nu$ and $\lambda $ as in the  statements of
Theorems \ref{eigensharp} and \ref{eigensharplambda},
respectively, there do exist a manifold of revolution $M$
fulfilling $\nu _M \approx \nu$ and a manifold of revolution $M$
fulfilling $\lambda _M \approx \lambda$. This is accomplished in
the following Proposition \ref{disisop}, dealing with $\lambda$,
and in Proposition \ref{disisoc}, dealing with $\nu$.

%
%

\begin{proposition}\label{disisop}
Let $n \geq 2$, and let $\lambda : [0, \infty ) \to [0, \infty )$
be
such that
\begin{equation}\label{n100}
\frac{\lambda (s)}{s^{1/n'}} \approx \hbox{ a non-decreasing
function near $0$.}
 \end{equation}
Then there exist $L \in (0, \infty ]$ and $\varphi : [0, L) \to
[0, \infty )$ as in the statement of Theorem \ref{disisopA} such
that:
\par \noindent (i) the $n$-dimensional manifold of revolution $M$
associated with $\varphi$ fulfills \eqref{l7bis} for some function
$\lambda _0$ such that
\begin{equation}\label{lambdabar}
\lambda _0 \approx \lambda \qquad \quad \hbox{near $0$};
\end{equation}
\par \noindent (ii)
the isoperimetric function $\lambda _M$of $M$ fulfills
\begin{equation}\label{eqlam}
\lambda _M \approx \lambda \qquad \quad \hbox{near $0$}.
\end{equation}
\par\noindent
Moreover, $L=\infty$ if and only if
\begin{equation}\label{l1}
\int _0 \frac{dr}{\lambda (r)} = \infty\,.
\end{equation}
\par\noindent
\end{proposition}
\begin{remark}\label{rem2}
{\rm If
\begin{equation}\label{noneiinftylambda}
\int _0 \frac r{\lambda (r)^2} dr = \infty\,,
\end{equation}
then \eqref{l1} holds. Indeed, if \eqref{l1} fails, namely if
\begin{equation}\label{**}
\int _0 \frac{dr}{\lambda (r)} < \infty\,,
\end{equation}
then  $\lim _{s\to 0}\frac {s}{\lambda
(s)}=0$, and
 this limit, combined with \eqref{**}, implies the convergence of the integral in
 \eqref{noneiinftylambda}.
 %
}
 \end{remark}
  \medskip
\par\noindent
 {\bf Proof of Proposition \ref{disisop}}.  Let $V$ be a positive number such that \eqref{n100} holds in $(0,
V)$; namely, there exists a non-decreasing function $\vartheta $
such that
$$\frac {\lambda (s)^{n'}}{s}
\approx \vartheta (s) \qquad \hbox{for $s \in (0, V)$.}$$ Thus,
the function $$\lambda _1 (s)= (s\vartheta (s))^{1/n'}$$ satisfies
$$\lambda _1 \approx \lambda  \qquad \hbox{in $(0, V)$,}$$
 and
\begin{equation}\label{lambda0}
\frac {\lambda _1 (s)^{n'}}{s} \qquad \hbox{is non-decreasing in
$(0, V)$.}
\end{equation}
 Assumption \eqref{lambda0} in turn ensures that, on
defining
$$\lambda _2 (s) = \bigg(\int _0^s \frac{\lambda _1(r)^{n'}}{r
} dr\bigg)^{1/n'} \qquad \hbox{for $s \in (0, V)$},$$ we have that
$\lambda _2 \in C^0(0 , V)$, $\lambda _2 ^{n'}$ is convex in $(0,
V)$, and $\lambda _2 \approx \lambda _1 \approx \lambda $ in $(0,
V)$. Similarly, on setting
$$\lambda _3 (s) = \bigg(\int _0^s \frac{\lambda _2(r)^{n'}}{r
} dr\bigg)^{1/n'} \qquad \hbox{for $s \in (0, V)$},$$
 we have that
$\lambda _3 \in C^1(0, V)$, $\lambda _3 ^{n'}$ is convex in $(0,
V)$, and $\lambda _3 \approx \lambda _2 \approx \lambda _1 \approx
\lambda$ in $(0, V)$.
\par\noindent
Thus, in what follows we may assume, on replacing if necessary
$\lambda $ by $\lambda _3$ near $0$, that $\lambda $ is a
non-decreasing function in $[0, L)$ such that $\lambda \in C^1(0,
V)$, $\lambda (0)=0$,
and
\begin{equation}\label{n100bis}
\lambda  ^{n'} \qquad \hbox{is convex in $(0, V)$.}
\end{equation}
Define
\begin{equation}\label{L}
L = 2\int _{0}^{V/2}\frac{dr}{\lambda (r)}\,,
\end{equation}
and note that $L= \infty $ if and only if \eqref{l1} is in force.
%
 Next, set
$$R = \begin{cases}
L/2 & \hbox{if $L< \infty$,}
\\ 1 & \hbox{if $L= \infty$.}
\end{cases}
$$
\par\noindent
 Let $N: [R, L ) \to
[0, V/2]$ be the function implicitly defined by
\begin{equation}\label{l2}
\int _{N(r)}^{V/2}\frac{dr}{\lambda (r)} = r-R \qquad \quad
\hbox{for $r \in [R, L)$\,.}
\end{equation}
Clearly, $N \in C^1(R,L)$ and $N$ decreases monotonically from
$V/2$ to $0$. Define $\varphi : [R, L ) \to [0, \infty )$ as
\begin{equation}\label{l3}
\varphi (r) = \bigg(\frac {\lambda (N(r))}{\omega
_{n-1}}\bigg)^{\frac 1{n-1}} \qquad \hbox{for $r \in (R, L)$\,,}
\end{equation}
and observe that $\varphi \in C^1(R,L)$. Since,
\begin{equation}\label{l4}
\lambda (N(r)) = -N'(r) \qquad \hbox{for $r\in [R, L)$}\,,
\end{equation}
and $N(L) =0$, one has that
\begin{equation}\label{l5}
\int _r^L \lambda (N(\rho)) d\rho = N(r) \qquad \hbox{for $r\in
[R, L)$}\,,
\end{equation}
whence
\begin{equation}\label{l6}
\lambda \bigg(\int _r^L \lambda (N(\rho)) d\rho\bigg) = \lambda
(N(r)) \qquad \hbox{for $r\in [R, L)$}\,,
\end{equation}
and finally, by \eqref{l3},
\begin{equation}\label{l7}
\lambda \bigg(\omega _{n-1}\int _r^L \varphi (\rho )^{n-1}d\rho
\bigg)= \omega _{n-1}\varphi (r )^{n-1} \qquad \hbox{for $r\in [R,
L)$}\,,
\end{equation}
namely \eqref{l7bis} with $\lambda _0$ replaced by $\lambda$.
\par\noindent
Now, observe that the function $\varphi $  is decreasing in $(R, L)$ and $\lim _{r \to L}\varphi (r) =0$. Furthermore, $ \varphi$ is convex
 owing to \eqref{n100bis}.
 Indeed, by \eqref{l4},
\begin{align}\label{l8}
\omega _{n-1}^{\frac 1{n-1}} \varphi '(r) = \big(\lambda
(N(r))^{\frac 1{n-1}}\big)'  & = \frac 1{n-1}\lambda
'(N(r))\lambda (N(r))^{\frac 1{n-1}-1}N'(r)\\ \nonumber & = -
\frac 1{n-1}\lambda '(N(r))\lambda (N(r))^{\frac 1{n-1}} \quad
\hbox{for $r \in (R, L)$.}
\end{align}
Thus, since $N(r)$ is decreasing, $\varphi '(r)$ is increasing if
and only if $- \lambda '(s)\lambda (s)^{\frac 1{n-1}}$ is
decreasing, namely if and only if $\lambda '(s)\lambda (s)^{\frac
1{n-1}}$ is increasing, and  this is  in turn equivalent to the
convexity of $\lambda (s)^{n'}$.
\par\noindent
Finally, let us continue $\varphi$ smoothly to the whole of $[0,
L)$ in such a way that \eqref{n0} and \eqref{n1'} are fulfilled,
and that $$\omega _{n-1} \int _0^R \varphi (r)^{n-1}\, dr = \omega
_{n-1} \int _R^L \varphi (r)^{n-1}\, dr = N(R) = \frac V2.$$
The resulting function $\varphi$ fulfils the
assumptions of Theorem \ref{disisopA}. Hence,
the conclusion  follows. \qed
\par\noindent
\smallskip
\begin{proposition}\label{disisoc}
Let $n \geq 2$, and let $\nu :[0, \infty ) \to [0, \infty )$ be a
function
such that
$\nu \in \Delta _2$ near $0$, and either $n \geq 3$ and
\begin{equation}\label{nu1}
\frac {\nu (s)}{s^{\frac{n-2}{n}}} \qquad \hbox{is equivalent to a non-decreasing function
near $0$,}
\end{equation}
or $n =2$ and there exists $\alpha >0$ such that
\begin{equation}\label{nu1bis}
\frac {\nu (s)}{s^\alpha } \qquad \hbox{is equivalent to a non-decreasing function
near $0$.}
\end{equation}
Then there exist $L \in (0, \infty ]$ and $\varphi : [0, L) \to
[0, \infty )$ as in the statement of Theorem \ref{disisopA}, such that:
\par\noindent
(i) the $n$-dimensional manifold of revolution $M$ built upon
$\varphi$ fulfills \eqref{l7bis} for some function $\lambda _0$
such that $\lambda _0 \approx \lambda _M$ near $0$;
\par\noindent
(ii)
\begin{equation}\label{nu0}
 \nu (s) \approx \nu _M (s) \approx \displaystyle{\frac{1}{\int _s^{\m2} \frac {dr}{\lambda
_M(r)^2}}}\qquad \hbox{near $0$}\,.
\end{equation}
\par\noindent
Moreover, $L =\infty $ if and only if
\begin{equation}\label{nu0'}
\int _0 \frac{dr}{\sqrt{r\nu (r)}} = \infty\,.
\end{equation}
\end{proposition}

\begin{remark}\label{rem3}
{\rm If
\begin{equation}\label{noneiinfty}
\int _0 \frac {dr}{\nu (r)}  = \infty\,,
\end{equation}
then \eqref{nu0'} holds. This is a consequence of the fact that
there exists an absolute constant $C$ such that  $$\bigg(\int _0^1
f(s)^2 ds\bigg)^{1/2} \leq C \int _0^1f(s) \frac {ds}{\sqrt s}$$
for every non-increasing function $f: (0, 1) \to [0, \infty )$.

%
%
}
 \end{remark}
{\bf Proof of Proposition \ref{disisoc}}. Let us assume that $ n \geq
3$, the case when $n=2$ being analogous. Let $V$ be a positive
number such that $\nu \in \Delta _2$ in $(0, V)$ and \eqref{nu1}
holds in $(0, V)$. An analogous argument as at the beginning of
the proof of Proposition \ref{disisop} tells us that on replacing $\nu$,
if necessary, by an equivalent function, we may assume that $\nu
\in C^1(0, V)$,
 $\nu '(s) >0$ for $s \in (0, V)$, and
 %
%
\begin{equation}\label{nu3}
s \nu ' (s) \approx \nu (s)   \qquad \hbox{for $s \in (0, V)$.}
\end{equation}
Define $\lambda : (0, V) \to (0, \infty )$ by
\begin{equation}\label{nu4}
\lambda (s) = \frac {\nu (s) }{\sqrt {\nu '(s)}}\qquad \hbox{for $s \in (0, V)$.}
\end{equation}
Given any $a \in (0, V)$, we thus have that
\begin{equation}\label{nu5}
\frac 1{\nu (s)} - \frac 1{\nu (a)} = \int _s^a \frac {dr}{\lambda
(r)^2} \qquad \hbox{for $s \in (0, V)$,}
\end{equation}
and hence there exists $\overline s$ such that
\begin{equation}\label{nu6}
\frac 1{2\nu (s)} \leq  \int _s^a \frac {dr}{\lambda (r)^2} \leq
\frac 1{\nu (s)} \qquad \hbox{if $0 < s <\overline s$.}
\end{equation}
Moreover,
\begin{equation}\label{n100ter}
\frac{\lambda (s)}{s^{1/n'}}  \quad\hbox{ is non-decreasing in
$(0, V)$.}
 \end{equation}
Indeed, by \eqref{nu3} and by the $\Delta _2$
condition for $\nu $ in $(0, V)$,
\begin{equation}\label{nu7}
\frac{\lambda (s)^{2}}{s^{\frac 2{n'}}} = \frac{\nu (s)^2} {\nu '
(s)s^{\frac 2{n'}}} \approx \frac{\nu (s)} {s^{\frac {n-2}{n}}}
\qquad \hbox{for  $s\in (0, V)$}\,.
\end{equation}
Owing to \eqref{nu7} and  \eqref{nu1}, the function $\lambda $
fulfills the assumptions of Proposition \ref{disisop}. Let
 $M$ be the $n$ dimensional manifold of revolution associated with $\lambda$
as in Proposition \ref{disisop}. By \eqref{lambdabar}, \eqref{eqlam},
\eqref{n50} and \eqref{n51},
\begin{align}\label{nu8}
\frac 1{\nu _M(s)} \approx
\int _{s}^{\m2} \frac {dr}{\lambda (r)^2} \qquad \hbox{near $0$.}
\end{align}
On the other hand, by \eqref{nu6},
\begin{equation}\label{nu1000}
 \int _{s}^{\m2} \frac {dr}{\lambda (r)^2} \approx \frac{1}{ \nu
(s)}\,\qquad \hbox{near $0$.}
\end{equation}
Equation \eqref{nu0} follows from \eqref{nu8} and \eqref{nu1000}.
%
\par
As for the assertion concerning \eqref{nu0'}, by Proposition
\ref{disisop} one has that $L=\infty$ if and only if the function $\lambda$ given by
\eqref{nu4} fulfils \eqref{l1}.
 One has that
 \begin{equation}\label{nu11}
 \int _0^V \frac {ds}{\lambda (s)} = \int _0^V \frac{\sqrt{\nu '
 (s)}}{\nu (s)} \geq C \int _0^V\frac{\sqrt {\nu (Cs)}}{\sqrt s \nu
 (s)}ds \geq C' \int _0 ^V\frac {ds}{\sqrt{s \nu (s)}}\,
 \end{equation}
 for suitable constants $C$ and $C'$, where the first inequality
 holds by \eqref{nu3} and the last one by the $\Delta _2$ condition
for $\nu$.  An analogous chain yields
$$\int _0^V \frac {ds}{\lambda (s)} \leq C \int _0 ^V\frac {ds}{\sqrt{s \nu (s)}}$$
for a suitable positive constant $C$.  Hence, equation
\eqref{nu0'} is equivalent to $L=\infty$.
%
\qed

\section{$L^q$ bounds for eigenfunctions}\label{sec5}
We begin  with the proof of Theorem \ref{eigencor} on $L^q(M)$
bounds for eigenfunctions.
\medskip
\par\noindent
{\bf Proof of Theorem \ref{eigencor}}. 
Given $s \in (0, \mo )$ and $h > 0$, choose the test function $v$
defined as
\begin{equation}\label{5'}
v (x)= \begin{cases}
0  & {\rm if}\,\,\, u (x)< u^\circ (s+h) \\
 u(x)-u^\circ (s+h) &
{\rm if}\,\,\, u^\circ (s+h) \leq u (x)\leq u^\circ (s)
\\
 u^\circ (s)-u^\circ (s+h) &
{\rm if}\,\,\, u^\circ (s) < u (x) \,,
\end{cases}
\end{equation}
in equation \eqref{weaksol}. Notice that $v \in W^{1,2}(M)$ by
standard results on truncations of Sobolev functions. We thus
obtain that
\begin{multline}\label{6}
\int _{\{u^\circ (s+h)<u<u^\circ (s)\}} |\nabla u|^2 \,d\hn (x)
 \\ =\gamma \int _{\{u^\circ (s+h) < u \leq u^\circ (s)\}} u(x) \big(u(x)-
u^\circ (s+h)\big)\,d\hn (x) \, + \, \big( u^\circ (s)-u^\circ
(s+h)\big) \gamma \int _{\{u
> u^\circ (s)\}} u (x)\,d\hn (x)\,.
\end{multline}
Consider the function $U: (0, \mo )\rightarrow [0, \infty)$ given
by
\begin{equation}\label{U}
U(s)=\int_{\{u \le u^\circ (s)\}} |\nabla u|^2\, d\hn (x) \qquad
\qquad \hbox{for $s\in (0, \mo )$.}
\end{equation}
By   \eqref{ac}, the function $u^\circ$ is locally absolutely
continuous (a.c., for short) in $(0, \mo )$.
The function
$$
(0, \infty) \ni t \mapsto \int_{\{u \le t\}}|\nabla u|^2\, d\hn
(x)
$$
is also locally a.c., inasmuch as, by the coarea formula,
\begin{equation}\label{coarea}
\int_{\{u \le t\}}|\nabla u|^2\, d\hn
(x)=\int_{-\infty}^t\int_{\{u= \tau\}} |\nabla u|d\hh (x) d\tau
\qquad \hbox{for } t\in \R.
\end{equation}
Consequently, $U$ is locally a.c., for it is the composition of
monotone locally a.c. functions, and
\begin{equation}\label{coareabis}
U'(s)=-{u^\circ}'(s)\int_{\{u= u^\circ(s)\}} |\nabla u|d\hh (x)
\quad \hbox{for a.e. $s\in(0, \mo )$}.
\end{equation}
Thus, dividing through by $h$ in \eqref{6},
and passing to the limit as $h \to 0^+$ yield
%
%
%
%
%
%
%
%
%
%
%
%
%
%
\begin{equation}\label{e6}
-{u^\circ}'(s)\int_{\{u= u^\circ(s)\}} |\nabla u|d\hh (x) = \gamma
(- {u^\circ} '(s))\int _{\{u>u^\circ (s)\}} u\, d\hn (x) \qquad \quad
\hbox{for a.e. $s \in (0, \mo )$.}
\end{equation}
On the other hand, it is easily verified via the definition of signed rearrangement that
\begin{equation}\label{phi1}
(- {u^\circ }'(s))\int _{\{u>u^\circ (s)\}}  u(x)\, dx = (-
{u^\circ }'(s))\int _0^s u ^\circ (r)\,dr \qquad \quad \hbox{for
a.e. $s \in (0, \mo )$.}
\end{equation}
Coupling \eqref{e6} and \eqref{phi1} tells us that
%

\begin{equation}\label{e8}
-{u^\circ} '(r)  = \frac{-{u^\circ} '(r)\,\gamma}
 {\int _{\{u=u^\circ (r) \}}
|\nabla u|\,d \hh (x)} \int _0^r u ^\circ (\varrho)\,d\varrho
\qquad \quad \hbox{for a.e. $r \in (0, \mo )$.}
%
%
%
\end{equation}
Let $0 < s \leq \varepsilon \leq \m2$. On integrating both sides
of \eqref{e8} over the interval $(s, \varepsilon )$, we obtain
that
\begin{equation}\label{e9}
u^{\circ}(s)  - \gamma \int _s^\varepsilon   \bigg(\int _0^{r}
 u ^\circ (\varrho)d\varrho\bigg)
 \big(-\psi _u (u^\circ (r))\big)' \,dr =
u^{\circ}(\varepsilon) \qquad \quad \hbox{for  $s \in (0,
\varepsilon )$,}
\end{equation}
where $\psi _u$ is the function defined as in \eqref{psi}.
 Define the operator $T$ as
\begin{equation}\label{e10}
Tf(s) = \int _s^\varepsilon   \bigg(\int _0^{r}
f(\varrho)d\varrho\bigg)\big(-\psi _u (u^\circ (r))\big)'\,dr
\qquad \quad \hbox{for  $s \in (0, \varepsilon )$,}
\end{equation}
for any integrable function $f$ on $(0, \varepsilon )$. Equation
\eqref{e9} thus reads
\begin{equation}\label{e11}
(I- \,\gamma T)(u^\circ ) = u^{\circ }(\varepsilon)\,.
\end{equation}
\par\noindent
We want now to show that, if  \eqref{1001} holds, then
 \begin{equation}\label{Tq}
T : L^q(0, \varepsilon ) \to L^q (0, \varepsilon ),
\end{equation}
and there exist an absolute constant $C$
 such that
\begin{equation}\label{normq}
\|T\|_{(L^q(0, \varepsilon ) \to L^q (0, \varepsilon ))} \leq C
\Theta (\varepsilon )\,,
\end{equation}
where $\Theta : (0, \m2 ] \to [0, \infty )$ is the function
defined as
\begin{equation}\label{theta}
\Theta (s ) = \sup _{r \in (0, s )} \frac{r}{\nu _M (r)}
 \qquad \quad \hbox{for $s \in (0, \m2 ]$.}
\end{equation}
Set
\begin{equation}\label{v}
v =  u- {\rm med}(u),
\end{equation}
 and observe that
\begin{equation}\label{medv}
{\rm med}(v) =0\,,
\end{equation}
$$v^\circ = u^\circ  - {\rm med}(u)\,,$$
and
\begin{equation}\label{e200}
\big(\psi _u (u^\circ (s))\big)'=\big(\psi _v (v^\circ (s))\big)'
\qquad \hbox{for $s \in (0, \mo )$.}
\end{equation}
%
 Moreover,
\begin{equation}\label{e11bis}
v^\circ (s) \geq 0 \qquad \quad \hbox{if $s \in (0, \m2 )$.}
\end{equation}
Given any  $q \in [2, \infty )$, $f \in L^q(0, \varepsilon )$ and
$0 <s \leq \varepsilon \leq \m2$, the following chain holds:
\begin{align}\label{e13}
|Tf(s)| & = \left|\int _s^\varepsilon   \bigg(\int _0^{r}
f(\varrho)d\varrho \bigg)\big(-\psi _u (u^\circ
(r))\big)'\,dr\right|
\\ \nonumber & =
\left|\int _s^\varepsilon   \bigg(\int _0^{r} f(\varrho)d\varrho
\bigg)\big(-\psi _v (v^\circ (r))\big)'\,dr\right|
\intertext{\qquad \qquad \qquad \qquad \qquad \qquad \qquad \qquad
\qquad \qquad \qquad \qquad \qquad \qquad \qquad (by
\eqref{e200})} &\leq  \int _s^\varepsilon \bigg(\int _0^{r}
|f(\varrho)|d\varrho\bigg)\frac{d}{dr}\Bigg( - \int _0^{v^\circ
(r)} \frac{d\tau}
 {\int _{\{v=\tau \}}
|\nabla v|\,d \hh (x)}\Bigg)\,dr \nonumber
\\
\nonumber & =  \bigg(\int _s^\varepsilon \frac{d}{dr}\Bigg( - \int
_0^{v^\circ (r)} \frac{d\tau}
 {\int _{\{v=\tau \}}
|\nabla v|\,d \hh (x)}\Bigg)\,dr \bigg)\int _0^s |f(\varrho )|d\varrho
\\ \nonumber & \quad + \int _s^\varepsilon \bigg(\int
_\varrho ^\varepsilon \frac{d}{dr}\Bigg(- \int _0^{v^\circ (r)}
\frac{d\tau}
 {\int _{\{v=\tau \}}
|\nabla v|\,d \hh (x)}\Bigg) dr \bigg) |f(\varrho )|d\varrho
\intertext{\qquad \qquad \qquad \qquad \qquad \qquad \qquad \qquad
\qquad \qquad \qquad \qquad \qquad \qquad \qquad (by Fubini's
theorem)} \nonumber & = \bigg(\int _{v^\circ
(\varepsilon)}^{v^\circ (s)} \frac{d\tau}
 {\int _{\{v=\tau \}}
|\nabla v|\,d \hh (x)}\bigg) \int _0^s  |f(\varrho )|d\varrho \\
\nonumber & \quad + \int _s^\varepsilon \int _{v^\circ
(\varepsilon)}^{v^\circ (\varrho )} \frac{d\tau}
 {\int _{\{v=\tau \}}
|\nabla v|\,d \hh (x)} |f(\varrho )|d\varrho
\\
\nonumber & \leq  \bigg(\int _{0}^{v^\circ (s)} \frac{d\tau}
 {\int _{\{v=\tau \}}
|\nabla v|\,d \hh (x)}\bigg)\int _0^s  |f(\varrho )|d\varrho \\
\nonumber & \quad +  \int _s^\varepsilon \int _{0}^{v^\circ
(\varrho )} \frac{d\tau}
 {\int _{\{v=\tau \}}
|\nabla v|\,d \hh (x)} |f(\varrho )|d\varrho \nonumber
\intertext{\qquad \qquad \qquad \qquad \qquad \qquad \qquad \qquad
\qquad \qquad \qquad \qquad \qquad \qquad \qquad ($v^\circ (\varepsilon ) \geq 0$ by
\eqref{e11bis})}\nonumber
\\ \nonumber &=
 \psi_v(v^\circ (s))\,\int _0^s  |f(\varrho )|d\varrho  +
\int _s^\varepsilon \psi_v(v^\circ (\varrho)) |f(\varrho
)|d\varrho
\\ \nonumber & \leq
 \frac 1{\nu _M(s)}\,\int _0^s  |f(\varrho )|d\varrho  +
 \int _s^\varepsilon \frac 1{\nu _M(\varrho )} |f(\varrho
)|d\varrho \intertext{\qquad \qquad \qquad \qquad \qquad \qquad
\qquad \qquad \qquad \qquad \qquad \qquad \qquad  (by \eqref{e0}
with $u$ replaced by $v$).}\nonumber
%
%
%
\end{align}
 Thus, if we show that
 there exist constants $C_1$ and $C_2$ such that
\begin{equation}\label{e202}
\bigg(\int _0^\varepsilon \bigg( \frac{1}{\nu _M(s)} \int _0^s
|f(r )|dr \bigg)^q ds \bigg)^{1/q} \leq C_1 \bigg( \int
_0^\varepsilon |f(s)|^q ds \bigg)^{1/q}
\end{equation}
and
\begin{equation}\label{e203}
\bigg(\int _0^\varepsilon \bigg(\int _s^\varepsilon \frac{1}{\nu
_M(r)} |f(r )|\,dr \bigg)^q ds \bigg)^{1/q} \leq C_2 \bigg( \int
_0^\varepsilon |f(s)|^q ds \bigg)^{1/q}
\end{equation}
for every $f\in L^q(0, \varepsilon )$, then we obtain that
\begin{equation}\label{e201}
\|T f\| _{L^q(0, \varepsilon )} \leq (C_1 +C_2) \|f\|_{L^q(0,
\varepsilon )}
\end{equation}
for every $f\in L^q(0, \varepsilon )$ By standard weighted Hardy
inequalities (see e.g. \cite[Section 1.3.1]{Ma2}), inequalities
\eqref{e202} and \eqref{e203} hold if and only if
\begin{equation}\label{e204}
\sup _{s \in  (0, \varepsilon )}s^{1/q'}\left\|1/\nu
_M\right\|_{L^q (s, \varepsilon )}<\infty
\end{equation}
and
\begin{equation}\label{e205}
\sup _{s \in  (0, \varepsilon )}s ^{1/q} \left\|{1}/{\nu _M}
\right\|_{L^{q'} (s, \varepsilon )} < \infty\,,
\end{equation}
respectively. Furthermore, the constants  $C_1$ and $C_2$ in
\eqref{e202} and \eqref{e203} are equivalent (up to absolute
multiplicative  constants) to the left-hand sides of \eqref{e204}
and \eqref{e205}, respectively.
\par\noindent
The left-hand sides of \eqref{e204} and \eqref{e205} agree if
$q=2$.  We claim that, if $q \in (2, \infty)$, then the left-hand
side of \eqref{e204} does not exceed the left-hand side of
\eqref{e205}, up to an absolute multiplicative constant. Actually,
since $\nu _M$ is non-decreasing,
\begin{align}\label{e6'}
\sup _{r\in (0,s)} r^{1/q}\left\|{1}/{\nu _M} \right\|_{L^{q'} (r,
s )} & \geq (s/2)^{1/q}\left\|{1}/{\nu _M} \right\|_{L^{q'} (s/2,
s )}
\\ & \geq   \frac{s}{2\nu _M(s)} \qquad \hbox{for $s \in (0, \m2 )$.}\nonumber
\end{align}
 Thus,
\begin{align}\label{e210}
s^{\frac{1}{q'}} \left\|{1}/{\nu _M} \right\|_{L^q (s ,
\varepsilon )} & \leq s^{\frac{1}{q'}}\bigg(\frac{1}{\nu _M(s)}
\bigg)^{1-\frac{q'}{q}} \left\|1/\nu _M\right\|_{L^{q'}(s ,
\varepsilon )}^{\frac{q'}{q}}  \\ \nonumber & =
\bigg(\bigg(\frac{s}{\nu _M(s)} \bigg)^{q-2}s^{1/q}\left\|1/\nu
_M\right\|_{L^{q'}(s , \varepsilon )}\bigg)^{\frac{1}{q-1}}
\\ & \leq 2^{\frac {q-2}{q-1}} \sup _{r \in  (0, \varepsilon )} r^{\frac{1}{q}} \left\|{1}/{\nu _M} \right\|_{L^{q'}(r ,
\varepsilon )} \nonumber \\ & \leq 2\sup _{r \in  (0, \varepsilon
)} r^{\frac{1}{q}} \left\|{1}/{\nu _M} \right\|_{L^{q'}(r ,
\varepsilon )}\, \qquad \qquad \hbox{for $s \in (0, \varepsilon
)$},\nonumber
\end{align}
where the last but one inequality holds owing to \eqref{e6'}.
Thus, our claim follows.  On the other hand,
\begin{align}\label{sup}
\sup _{s \in  (0, \varepsilon )} & s ^{1/q} \left\|{1}/{\nu _M}
\right\|_{L^{q'} (s, \varepsilon )} = \sup _{s \in  (0,
\varepsilon )} s ^{1/q} \bigg(\int _s^\varepsilon
\bigg(\frac{r}{\nu _M (r)}\bigg)^{q'} \frac{dr}{r^{q'}}
\bigg)^{1/q'} \\ \nonumber &\leq \bigg(\sup _{s \in  (0,
\varepsilon )} \frac{s}{\nu _M (s)} \bigg) \sup _{s \in  (0,
\varepsilon )} s ^{1/q} \bigg(\int _s^\varepsilon
\frac{dr}{r^{q'}} \bigg)^{1/q'} \leq \frac 1{(q' -1)^{1/q'}}
\Theta (\varepsilon )\,.
\end{align}
Hence,  \eqref{normq} follows as well.
 \par\noindent
 On denoting by  $(I-\, \gamma  T)_q$ the restriction of
 $I-\,\gamma
T$ to $L^q (0, \varepsilon )$,  we deduce
via a classical result of functional analysis
 that the operator
\begin{equation}\label{e211}
(I- \gamma T)_q : L^q (0, \varepsilon ) \to L^q (0, \varepsilon )
\end{equation}
is invertible, with a bounded inverse, provided that $\varepsilon$
is so small that
\begin{equation}\label{cgamma}
C\gamma \Theta (\varepsilon ) < 1\,,
\end{equation}
where $C$ is the constant appearing in \eqref{normq}.
 Moreover,
\begin{equation}\label{e212}
\|(I-\gamma T)_q ^{-1}\| \leq \frac{1}{1-C\gamma\Theta
(\varepsilon)}\,.
\end{equation}
\par
The next step consists in showing that there exists an absolute
constant $C'$ such that also
 the restriction $(I- \, \gamma T)_{2}$ of $I- \, \gamma T$ to $L^{2} (0, \varepsilon )$ is invertible,
with a bounded inverse, provided that
\begin{equation}\label{c'gamma}
C'\gamma \Theta (\varepsilon ) < 1\,.
\end{equation}
An analogous chain as in \eqref{e210} tells us that
\begin{equation}\label{e210bis}
s^{1/2}\left\| {1}/{\nu _M} \right\| _{L^2(s, \varepsilon )} \leq
2 \sup _{r \in (0, \varepsilon )} r^{1/q}\left\| {1}/{\nu _M}
\right\| _{L^{q'}(r, \varepsilon )}
 \qquad \hbox{for $s
\in (0, \varepsilon )$}.
\end{equation}
 Consequently,
\begin{equation}\label{e213}
\sup _{s\in (0,\varepsilon )} s^{1/2}\left\| {1}/{\nu _M} \right\|
_{L^{2}(s, \varepsilon )}  \leq 2\Theta (\varepsilon )\,.
\end{equation}
Hence, the invertibility of $(I- \, \gamma T)_{2}$ under
\eqref{c'gamma} follows via the same argument as above. Owing to
assumption \eqref{1001}, both \eqref{cgamma} and \eqref{c'gamma}
hold provided that $\varepsilon$ is sufficiently small. In
particular, one can choose
\begin{equation}\label{eps}
 \varepsilon = \Theta ^{-1}\big(1/(2 \gamma C'')\big) \,,
 \end{equation}
 where $\Theta ^{-1}$ is the generalized  left-continuous inverse
 of $\Theta$, and $C'' = \max\{C,
 C'\}$.
%
\par\noindent
We have that $u^\circ \in  L^{2} (0,
\varepsilon )$, for $u \in L^2(M)$. Thus, since
 the constant function $u^\circ (\varepsilon )$ trivially belongs to $L^q
(0, \varepsilon )\subset L^{2} (0, \varepsilon )$, from
\eqref{e11} we deduce that
\begin{equation}\label{e11sette}
u^\circ = (I- \gamma T)_{2}^{-1} (u^\circ (\varepsilon )) = (I-\,
\gamma  T)_q^{-1} (u^\circ (\varepsilon ))\,.
\end{equation}
Hence, $u^\circ \in  L^q (0, \varepsilon )$, and, by \eqref{e212}
and \eqref{eps},
\begin{equation}\label{e12}
\|u^\circ \|_{L^q (0, \varepsilon )} \leq \frac{\|u^{\circ
}(\varepsilon)\|_{L^q(0,\varepsilon )}}{1-C \gamma \Theta
(\varepsilon )} = \frac{\varepsilon ^{1/q}|u^{\circ
}(\varepsilon)|}{1-C \gamma \Theta (\varepsilon )} \leq 2
\varepsilon ^{1/q}|u^{\circ }(\varepsilon)|\,.
\end{equation}
Since $\varepsilon \leq \m2$, one can easily verify that
\begin{equation}\label{n1} \|u \|_{L^{2} (M)} \geq
 \varepsilon
^{1/{2}}|u^\circ (\varepsilon )|\,.
\end{equation}
From \eqref{e12} and \eqref{n1}
 one has that
\begin{equation}\label{n2bis}
\|u^\circ \|_{L^q (0, \varepsilon )} \leq \frac{2\|u \|_{L^{2}
(M)}}{\varepsilon ^{\frac 1{2}- \frac 1q}}
\,.
\end{equation}
Next, observe that
\begin{equation}\label{med0}
|{\rm med }(u)| \leq (2/{\mo })^{1/2} \|u\|_{L^2(\Omega )}.
\end{equation}
By \eqref{n1} and \eqref{med0}, there exists a constant $C=C(\mo
)$ such that
\begin{align}\label{med1}
\|u^\circ - {\rm med}(u)\|_{L^q(\varepsilon , \m2 )} & \leq
(u^{\circ} (\varepsilon ) - {\rm med }(u))^{\frac
{q-2}{q}}\|u^\circ - {\rm med}(u)\|_{L^2(\varepsilon , \m2
)}^{\frac 2q} \\ \nonumber & \leq C (\varepsilon ^{\frac 1{2}-
\frac 1q} + 1 )\|u \|_{L^{2} (M)}.
\end{align}
From \eqref{n1}, \eqref{med0} and \eqref{med1} we deduce that
\begin{align}\label{med2}
\|u^\circ \|_{L^q(0, \m2 )} & \leq \|u^\circ \|_{L^q(0,
\varepsilon )} + \|u^\circ - {\rm med}(u)\|_{L^q(\varepsilon , \m2
)} + \|{\rm med}(u) \|_{L^q(\varepsilon , \m2 )} \\ \nonumber &
\leq (\varepsilon ^{\frac 1{2}- \frac 1q} + 1 )\|u \|_{L^{2} (M)}
\end{align}
for some constant $C=C(\mo )$. Hence, there exists a constant
$C=C(\mo )$ such that
\begin{align}\label{q7bis}
\|u^\circ\|_{L^q (0, \m2 )}  & \leq \frac{C \|u \|_{L^{2} (M)}}{
\big(\Theta ^{-1}\big(1/( \gamma C)\big)\big)^{\frac 1{2}- \frac
1q} }\,.
\end{align}

%
%
%
%
%
%
%
%
%
A combination of  \eqref{med2}  with  an analogous estimate for
 $\|u^\circ\|_{L^q (\m2 ,
\mo )}$ obtained via the same argument applied to $-u$,  yields
\eqref{ei5}, since $(-u)^\circ (s) = -u^\circ (\hn (M) -s)$ for $s
\in (0, \hn (M))$.
\qed
\medskip
\par
An inspection of the proof of Theorem \ref{eigencor} reveals that,
in fact,  the following estimate for the constant appearing in
\eqref{ei5} holds.

\begin{proposition}\label{eigencorconstant}
Define the function $\Theta : (0, \m2 ] \to [0, \infty )$ as
\begin{equation*}
\Theta (s ) = \sup _{r \in (0, s )} \frac{r}{\nu _M (r)}
 \qquad \quad \hbox{for $s \in (0, \m2 ]$.}
\end{equation*}
Then inequality \eqref{ei5} holds with
$$C(\nu _M , q, \gamma ) = \frac{C_1}{\big(\Theta ^{-1}(C_2/\gamma )\big)^{\frac 12
-\frac 1q}}\,,$$ where $C_1=C_1(q, \mo )$ and $C_2=C_2(q, \mo )$
are suitable constants, and $\Theta ^{-1}$ is the generalized
left-continuous inverse of $\Theta$.
\end{proposition}

\begin{example}\label{example1}
{\rm Assume that there exists $\beta \in [(n-2)/n,1)$ such that
the manifold $M$ fulfils $\nu _M (s) \geq C s ^\beta$ for some
positive constant $C$ and for small $s$. Then \eqref{1001} holds,
and, by Proposition \ref{eigencorconstant}, for every $q \in (2,
\infty )$ there exists an  constant $C=C(q, \mo )$  such that
\begin{equation*}
\|u \|_{L^q (M )} \leq C \gamma ^{\frac{q-2}{2q(1-\beta )}}
 \|u \|_{L^{2} (M )}\,
\end{equation*}
for every eigenfunction $u$ of the Laplacian on $M$ associated with the eigenvalue $\gamma$.
}
\end{example}

\bigskip
\par\noindent
We now prove Theorem \ref{Lqsharp}.
\medskip
\par\noindent
{\bf Proof of Theorem \ref{Lqsharp}} Given $q>2$ and $ n\geq 2$,
we shall construct an $n$-dimensional  manifold of revolution $M$
as in Theorem \ref{disisopA}, fulfilling \eqref{borderline} and
such that the Laplacian on $M$ has an eigenfunction $u \notin
L^q(M)$. \par\noindent In the light of the discussion preceding
Theorem \ref{disisopA}, in order to exhibit such  eigenfunction it
suffices to produce a positive number $\gamma $ and a smooth
function $p : \R \to (0, \infty )$ such that
$$\int _{-\infty}^s \sqrt{p(\varrho )} d\varrho < \infty \qquad \quad \hbox{for $s \in \R$,}$$
and
$$\int _{\R} \sqrt{p(\varrho )} d\varrho = \infty \,,$$
 and a function $v : \R \to \R$ fulfilling \eqref{n109} and
\eqref{n111}, and such that
\begin{equation}\label{q1}
\int _\R v(s)^q p(s) ds = \infty\,.
\end{equation}
The function $\varphi$ in the definition of $M$ is recovered from
$p$ by
\begin{equation}\label{q2}
\varphi (r) = p(F^{-1}(r))^{\frac 1{2(n-1)}} \qquad \hbox{for
$r>0$},
\end{equation}
and $\varphi (0)=0$, where $F: \R \to [0, \infty )$ is given by
\begin{equation}\label{q3}
F(s) = \int _{-\infty}^s \sqrt{p(\varrho )} d\varrho  \qquad
\hbox{for $s \in \R$}.
\end{equation}

\par
We define the function $p$ piecewise as follows. Let $s_1 \leq -1
\leq 1 \leq s_2$ to be fixed later, and set
\begin{equation}\label{q4}
p(s) = \frac 1{s^2} \quad \hbox{if $s \geq s_2$}.
\end{equation}
Let $$0 <\gamma < \frac 14,$$ and
$$\alpha = \frac {1-\sqrt{1-4\gamma}}{2}.$$
With this choice of $\alpha$, the function
\begin{equation}\label{q6}
v(s) = s ^\alpha
\end{equation}
solves equation \eqref{n109} in $[s_2, \infty )$. On the other
hand, if   $p$ is defined in  $(- \infty, s_1]$ by
\begin{equation}\label{q5}
p(s)=
\begin{cases} \frac{4e^{2s}}{\gamma (1- e^{2s})} & \qquad \hbox{if
$n=2$,}
\\
\frac{ (-s)^{\frac{2n-2}{2-n}}}{(n-2)^{\frac{2n-2}{n-2}}-
\frac{\gamma (n-2)^2}{2n} (-s)^{\frac 2{2-n}}} & \qquad \hbox{if
$n>2$,}
\end{cases}
\end{equation}
then the function given by
\begin{equation}\label{q7}
 v(s) =\begin{cases} 1- \,e^{2s} & \qquad
\hbox{if $n=2$,}
\\
(n-2)^{\frac{2n-2}{n-2}} - \frac{\gamma (n-2)^2}{2n}(-s)^{\frac
2{2-n}} & \qquad \hbox{if $n>2$,}
\end{cases}
\end{equation}
solves equation \eqref{n109} in $(- \infty, s_1]$.
%
%
Next, given $\beta >0$ and neighborhoods $I_{-1 }$ and
 $I_{1 }$ of $-1$ and $1$, respectively, let $p$ be
defined in $I_{1 }\cup I_1$ as
$$p(s) =
\begin{cases}
\frac 6{\gamma (\beta - (s- 1)^2)} & \hbox{for $s \in I_{1 }$,}
\\ \frac 6{\gamma (\beta - (s+ 1)^2)} & \hbox{for $s \in I_{-1 }$}.
\end{cases}
$$
Then the function $v$ given by
$$
v(s) =
\begin{cases}
\beta (s- 1) - (s- 1)^3 & \hbox{for $s \in I_{1 }$,}\\
- \beta (s+ 1) + (s+ 1)^3 &  \hbox{for $s \in I_{-1 }$,}
\end{cases}
$$
 is a solution to \eqref{n109} in $I_{1} \cup I_1$.
%
Moreover, $v$ is convex in a left neighborhood of $1$ and in a
right  neighborhood of $-1$, whereas it is concave in a right
neighborhood of $1$ and in a left  neighborhood of $-1$. It is
easily seen that, if $\beta$ is sufficiently large, $s_2$ and
$-s_1$ are sufficiently large depending on $\beta$,
 and $I_{1}$ and $I_{-1}$ have a sufficiently
small radius, then $v$ can be continued to the whole of $\R$ in
such a way that:
 $$\hbox{$v \in C^2(\R)$;}$$
$$\hbox{$ v''  \leq - C$ and   $v \geq C$ in $\R \setminus (I_{-1} \cup
(-1 , 1) \cup I_{1 })$, for some positive constant C;}$$
 $$\hbox{ $v''  \geq  C$ and   $v \leq - C$ in $(-1, 1) \setminus  (I_{-1} \cup I_{1
  })$, for some positive constant C.}$$
%
%
%
\par\noindent
Thus, $p$ can be continued to the whole of $\R$ as a positive
function in $C^2(\R )$ in such a way that $v$ is a solution to
equation \eqref{n109} in $\R$.
Also, the function $v$  fulfils \eqref{n111} for every $\gamma \in
(0, \frac 14)$, and  \eqref{q1} provided that $\gamma$ is
sufficiently close to $\frac 14$.
\par
The manifold of revolution $M$ built upon the function $\varphi$
given by \eqref{q2} satisfies the assumtpions of Theorem
\ref{disisopA}. Indeed, $\varphi (r) >0$ if $r>0$, and \eqref{n1'}
holds, as a consequence of the fact that
%
\begin{equation}\label{limite}
\lim _{s \to - \infty } \frac{p'(s)}{ p(s)^{\frac{3n -4}{2n-2}}} =
2(n-1)\,.
\end{equation}
%
%
Assumptions (i)--(iii) of Theorem \ref{disisopA} are also
fulfilled, since there exists $b>0$ such that
\begin{align}\label{exp}
\varphi (r) = b e^{-\frac r{n-1}} \qquad \quad \hbox{ for large
$r$.} \end{align}
 Finally, if $\lambda _0$ denotes the function given by \eqref{l7bis}, then  we
 obtain via \eqref{exp} that
 \begin{align}\label{b}
\lim _{s \to 0}\frac{s}{ \lambda _0(s )}
 \approx \lim _{r\to \infty } \frac 1{\varphi (r)^{n-1}}\int _r^\infty
\varphi (\rho )^{n-1}d\rho =1.
\end{align}
Hence, by \eqref{n51},
\begin{align*}
\lim _{s \to 0} \frac{s}{\nu _M (s)} & \approx \lim _{s \to 0}
s\int _s^{\m2}\frac {d\varrho}{ \lambda _0(\varrho )^2} \approx
\lim _{s \to 0}\frac{s^2}{ \lambda _0(s )^2} \approx 1.
\end{align*}
Hence \eqref{borderline} follows. \qed
\smallskip
\par
We conclude this section by sketching  the proofs of the results
dealing with $L^q(M)$ estimates in terms of $\lambda _M$.
%
%
%

%
%

%
%
\smallskip
\par\noindent
{\bf Proof of Theorem \ref{Lqlambdabis}} By \eqref{1001bis}, given
$\varepsilon >0$, there exists $s_\varepsilon$ such that
$\frac{s^2}{\lambda _M (s)^2} <\varepsilon$ if $0<s<s
_\varepsilon$. By inequality \eqref{lamnu}, if $0<s<s
_\varepsilon$, then
\begin{align}\label{a}
\frac {s}{\nu _M (s)} \leq s \int _{s}^{\m2 } \frac{dr}{\lambda _M
(r)^2}
 & \leq \varepsilon
s \int _{s}^{s _\varepsilon }\frac {dr}{r^2}+ s \int _{s
_\varepsilon}^{\m2 } \frac{dr}{\lambda _M (r)^2} \\ \nonumber &
\leq \varepsilon +s \int _{s _\varepsilon}^{\m2 }
\frac{dr}{\lambda _M (r)^2}.
\end{align}
The rightmost side of \eqref{a} tends to $\varepsilon$ as $s \to
0$. Hence \eqref{1001} follows, owing to the arbitrariness of
$\varepsilon$. Inequality \eqref{ei5lambda} in thus a consequence
of Theorem \ref{eigencor}. \qed
\smallskip
\par\noindent
{\bf Proof of Theorem \ref{Lqsharplambda}} The proof is  the same
as that of Theorem \ref{Lqsharp}. One has just to notice that
\eqref{b} and \eqref{n50} imply \eqref{borderlinelambda}. \qed

\section{Boundedness of eigenfunctions}\label{sec4}

Our main task in this section is to prove Theorem
\ref{eigencorinf}, which provides a condition on $\nu _M$ for the
boundedness of the eigenvalues of the Laplacian, and  Theorem
\ref{eigensharp}, showing the sharpness of such  condition. The
proofs of the parallel results of Theorems \ref{boundlambda} and
\ref{eigensharplambda} involving $\lambda _M$ are sketched at the
end of the section.
\smallskip
\par\noindent
 {\bf Proof of  Theorem \ref{eigencorinf}} We start as in the
proof of Theorem  \ref{eigencor},  define the operator $T$ as in
\eqref{e10}, and, for $\varepsilon \in (0, \m2 )$, we write again
equation \eqref{e9} as
\begin{equation}\label{eqinf}
(I-\gamma \,T)(u^\circ ) = u^{\circ }(\varepsilon)\,.
\end{equation}
\par\noindent
We claim that, if  \eqref{eiinftylambda} is
satisfied,  then
 \begin{equation}\label{Tinf}
T : L^\infty (0, \varepsilon ) \to L^\infty (0, \varepsilon ),
\end{equation}
and
\begin{equation}\label{norminf}
\|T\|_{(L^\infty (0, \varepsilon ) \to L^\infty (0, \varepsilon
))} \leq \int _0^\varepsilon \frac {dr}{\nu _M (r)}\,.
\end{equation}
To verify our claim,  define $v$ as in \eqref{v}, recall
\eqref{e200}, and note that
\begin{align}\label{infe13}
\|Tf\|_{L^\infty (0, \varepsilon )} & \leq  \|f\|_{L^\infty (0,
\varepsilon )}\int _0^\varepsilon \bigg(\int _0^{r}
d\varrho\bigg)\frac{d}{dr}\Bigg( -\int _0^{v^\circ (r)} \frac{d\tau}
 {\int _{\{v=\tau \}}
|\nabla v|\,d \hh (x)}\Bigg)\,dr \\ \nonumber & =  \|f\|_{L^\infty
(0, \varepsilon )} \int _0^\varepsilon \int _\varrho ^\varepsilon
\frac{d}{dr}\Bigg( -\int _0^{v^\circ (r)} \frac{d\tau}
 {\int _{\{v=\tau \}}
|\nabla v|\,d \hh (x)}\Bigg)\,dr\, d\varrho \\ \nonumber &=
\|f\|_{L^\infty (0, \varepsilon )} \int _0^\varepsilon \int
_{v^\circ (\varepsilon )}^{v^\circ (\varrho )} \frac{d\tau}
 {\int _{\{v=\tau \}}
|\nabla v|\,d \hh (x)}d\varrho \\ \nonumber & \leq \|f\|_{L^\infty
(0, \varepsilon )} \int _0^\varepsilon \int _{0}^{v^\circ (\varrho
)} \frac{d\tau}
 {\int _{\{v=\tau \}}
|\nabla v|\,d \hh (x)}d\varrho
\\ \nonumber & =
\|f\|_{L^\infty (0, \varepsilon )}  \int _0^\varepsilon \psi _v
(v^\circ (\varrho )) \, d\varrho
 \\ \nonumber & \leq \|f\|_{L^\infty
(0, \varepsilon )} \int _0^\varepsilon \frac{d\varrho}{\nu
_M(\varrho)}.
%
\end{align}
Observe that we have made use of the inequality $v(\varepsilon )
\geq 0$, due to \eqref{e11bis}, in the last but one inequality,
and of \eqref{e0} (with $u$ replaced by $v$) in the last
inequality.
Now, define the function $\Xi : (0, \m2 ] \to [0, \infty )$ as
\begin{equation}\label{xi}
\Xi (s ) =  \int _0^s \frac {dr}{\nu _M (r)}
 \qquad \quad \hbox{for $s \in (0, \m2 ]$.}
\end{equation}
If  $\varepsilon$ is sufficiently  small, in particular
$$\varepsilon \leq \Xi ^{-1}\big(1/(2\gamma )\big),$$ where $ \Xi
^{-1}$ is the generalized left-continuous inverse of $ \Xi $, then
%
%
 we deduce from \eqref{norminf} that the
restriction $(I-\gamma \,T)_\infty$ of $I-\gamma \, T$ to
$L^\infty (0, \varepsilon )$,
\begin{equation}\label{e11'}
(I- \gamma \,T)_\infty : L^\infty (0, \varepsilon ) \to L^\infty
(0, \varepsilon )
\end{equation}
is invertible, with a bounded inverse, and
\begin{equation}\label{e11''}
\|(I- \gamma \,T)_\infty ^{-1}\| \leq \frac{1}{1-\gamma \int
_0^\varepsilon \frac{d\varrho}{\nu _M(\varrho)}}\leq 2\,.
\end{equation}
Since $\Xi (\varepsilon ) \geq \Theta (\varepsilon )$, where
$\Theta $ is the function defined by \eqref{theta}, we infer from
\eqref{e213} and from the proof of Theorem \ref{eigencor} that
there exists an absolute constant $C$ such that the restriction
%
\begin{equation}\label{e11sei'}
(I- \gamma \,T)_2 : L^2 (0, \varepsilon ) \to L^2 (0, \varepsilon
)
\end{equation}
of $I-\gamma T$ to $L^2 (0, \varepsilon )$ is also invertible,
with a bounded inverse, provided that $C\Xi (\varepsilon ) <1$.
Set $C' = \max\{2, C\}$, and choose
$$\varepsilon = \Xi^{-1} \big(1/(\gamma C')\big).$$
%
%
%
%
%
\par\noindent
Since $u^\circ \in  L^2 (0, \varepsilon )$ and $u^\circ
(\varepsilon ) \in L^\infty (0, \varepsilon )\subset L^2 (0,
\varepsilon )$, from \eqref{eqinf} we deduce that
\begin{equation}\label{e11sette'}
u^\circ = (I- \gamma \,T)_2^{-1} (u^\circ (\varepsilon )) = (I-
\gamma T)_\infty^{-1} (u^\circ (\varepsilon ))\,.
\end{equation}
Hence, $u^\circ \in  L^\infty (0, \varepsilon )$, and
\begin{equation}\label{e12'}
\|u^\circ \|_{L^\infty (0, \varepsilon )} \leq \frac{|u^{\circ
}(\varepsilon)|}{1-\gamma  \int _0^\varepsilon \frac{d\varrho}{\nu
_M(\varrho)}}\leq 2|u^{\circ }(\varepsilon)|\,.
\end{equation}
Since
\begin{equation}\label{m1} \|u \|_{L^2 (M)}  \geq \varepsilon ^{1/2}|u^\circ
(\varepsilon )|=  \big(\Xi^{-1} \big(1/(\gamma
C')\big)\big)^{1/2}|u^\circ (\varepsilon )|\,,\end{equation}
 we have that
\begin{equation}\label{m2}
u^\circ (0) \leq \|u^\circ \|_{L^\infty (0, \varepsilon )} \leq
\frac{2\|u \|_{L^2 (M)}}{\big(\Xi^{-1} \big(1/(\gamma
C')\big)\big)^{1/2}}\,.\end{equation}
 The same argument,
applied to $-u$,
yields the same estimate for $-u^\circ (\mo )$.
 Since
$$
\|u \|_{L^\infty (M)} = \max \{u^\circ (0), - u^\circ (\mo)\},$$
inequality \eqref{ei5'} follows.
\qed

The following estimate for the constant  in \eqref{ei5'} is
provided in the proof of Theorem \ref{eigencorinf}.

\begin{proposition}\label{eigencorinfconstant}
Assume that \eqref{eiinftylambda} is in force.
Define the function $\Xi : (0, \m2 ] \to [0, \infty )$ as
\begin{equation*}
\Xi (s ) =  \int _0^s \frac {dr}{\nu _M (r)}
 \qquad \quad \hbox{for $s \in (0, \m2 ]$.}
\end{equation*}
Then inequality \eqref{ei5'} holds with
$$C(\nu _M ,  \gamma ) =
 \frac{C_1}{\big(\Xi ^{-1}(C_2/\gamma )\big)^{
\frac 12}}\,,
$$
 where $C_1$ and $C_2$ are suitable absolute
constants, and $\Xi ^{-1}$ is the generalized left-continuous
inverse of $\Xi$.
\end{proposition}

\begin{example}\label{example2}
{\rm Assume that there exists $\beta \in [(n-2)/n,1)$ such that
the manifold $M$ fulfils $\nu _M (s) \geq C s ^\beta$  for some
positive constant $C$ and for small $s$. Then
\eqref{eiinftylambda} holds, and, by Proposition
\ref{eigencorinfconstant}, there exists an absolute constant $C$
such that
\begin{equation*}
\|u \|_{L^\infty (M )} \leq C \gamma ^{\frac{1}{2(1-\beta)}}
 \|u \|_{L^{2} (M )}\,
\end{equation*}
for every eigenfunction $u$ of the Laplacian on $M$ associated with the eigenvalue $\gamma$.
}
\end{example}

%
%
%
\bigskip
\par\noindent
Next, we give a proof of Theorem \ref{eigensharp}.
\smallskip
\par\noindent
 {\bf Proof of  Theorem \ref{eigensharp}}
By Proposition \ref{disisoc},   if $\nu$ is as in the statement,
%
%
%
then there exists a function $\varphi$ such that the associated
$n$-dimensional manifold of revolution $M$  (as in Section
\ref{rev}) fulfils \eqref{1004}, and hence
\begin{equation}\label{sendai1}
\lim _{s \to 0}\frac s{\nu _M (s)} = \lim _{s \to 0}\frac s{\nu
(s)}=0
\end{equation}
 and
 \begin{equation}\label{sendai2}
\int _0 \frac {ds}{\nu _M (s)} = \int _0 \frac {ds}{\nu  (s)}=
\infty\,.
\end{equation}
In particular, the latter equation entails, via Remark \ref{rem3},
that \eqref{nu0'} holds, and hence that $L=\infty$ in Proposition
\ref{disisoc}. Thus, $\varphi : [0, \infty ) \to [0, \infty )$.
%
%
%
%
%
Now, recall that the function $\varphi $ given by Proposition
\ref{disisoc} is defined in such a way that \eqref{n1'} holds.
Hence,
\begin{equation}\label{n102'}
\int _0^1 \frac{dr}{\varphi (r )^{n-1}} = \infty\,,
\end{equation}
and
\begin{equation}\label{n102''}
\lim _{r\to 0} \bigg(\int _r^1 \frac{d\rho}{\varphi (\rho
)^{n-1}}\bigg)\bigg( \int _0^r\varphi (\rho )^{n-1} d\rho\bigg)
=0\,.
\end{equation}
Moreover,
\begin{equation}\label{n102'''}
\int _1^\infty \frac{dr}{\varphi (r )^{n-1}} = \infty\,,
\end{equation}
since $\lim _{r \to \infty } \varphi (r) = 0$ by property (i) of
Theorem \ref{disisopA}.
%
%
%
\par\noindent
Owing to Corollary \ref{manifoldsofrevolution},  condition
\eqref{sendai1} is equivalent to
\begin{equation}\label{limite0}
\lim _{r\to \infty } \bigg(\int _1^r\frac {d\varrho}{\varphi
(\varrho)^{n-1}}\bigg)\bigg(\int _r^\infty \varphi (\varrho
)^{n-1}d\varrho\bigg) =0\,,
\end{equation}
%
and condition \eqref{sendai2} is equivalent to
\begin{equation}\label{n104}
\int ^\infty \bigg(\frac 1{\varphi (r)^{n-1}}\int _r^\infty
\varphi (\rho )^{n-1}d\rho\,\bigg)\, dr = \infty\,.
\end{equation}
By the discussion preceding Theorem \ref{disisopA}, the conclusion
will follow if we exhibit   a number $\gamma
>0$ and an unbounded solution $v : \R \to \R$ to equation
\eqref{n109} fulfilling \eqref{n111}.
%
%
%
Note that $s_0 = \infty $ in \eqref{*} owing to \eqref{n102'''}.
%
Conditions \eqref{n102''}   and \eqref{limite0} are equvalent to
\begin{equation}\label{n112bis}
\lim _{s \to -\infty} s \int _{-\infty}^s p(t) dt =0\,,
\end{equation}
and
\begin{equation}\label{n112}
\lim _{s \to \infty} s \int _s ^{\infty} p(t) dt =0\,,
\end{equation}
respectively.
Condition \eqref{n104} amounts to
\begin{equation}\label{n112,'}
 \int  ^{\infty }s\,p(s)\, ds =\infty \,.
\end{equation}
Assumptions \eqref{n112bis} and \eqref{n112} ensure that  the
 the embedding
\begin{equation}\label{n113}
W^{1,2}(\R ) \to L^2(\R , p(s)ds)
\end{equation}
is compact -- see e.g. \cite{KK}. Here, $ L^2(\R , p(s)ds)$
denotes the space $L^2$ on $\R$ equipped with the measure
$p(s)ds$.
Consider the functional
\begin{equation}\label{n114}
J(v) = \frac {\int _\R \big(\frac{dv}{ds}\big)^2ds}{\int _\R
v^2\,p(s) ds}\,.
\end{equation}
 By
 the
compactness of embedding \eqref{n113}, there exists
$\min  J(v)$ among
all (non trivial) functions $v\in W^{1,2}(\R )$ such that $\int
_\R v\,p(s) ds=0$. Moreover, the minimizer $v$ is a solution to
equation \eqref{n109} with $\gamma = \min J$.
\par
By Hille's theorem \cite{Hille}, condition \eqref{n112} entails
that equation \eqref{n109} is nonoscillatory at infinity, and
hence that every solution has constant sign at infinity. Thus, we
may assume that $v(s)>0$ for large $s$. Consequently,
$$\frac{d^2v}{ds^2} < 0 \qquad \hbox{for large $s$,}$$
and hence $v$ is concave near $\infty$. Now, assume by
contradiction that $v$ is bounded. Then there exists $\lim _{s\to
\infty} v(s)$, and, on denoting by $v(\infty )$ this limit, one
has that
 $v(\infty ) \in (0, \infty )\,.$
Moreover,
$$
\lim _{s\to \infty} \frac{dv}{ds} =0\,.$$ Integration of
\eqref{n109} and the last equation yield
$$\frac{dv}{ds} = \gamma \int _s^\infty p(t) v(t) dt  \qquad \hbox{for large $s$.}
$$
Hence, there exists $\widehat s$ such that
$$ \frac{dv}{ds} \geq  \gamma \frac{v (\infty )}{2}\int _s^\infty p(t)  dt
\qquad \quad \hbox{ for $s \geq \widehat s$.}$$ By a further
integration, we obtain that
$$v(\infty ) - v(\widehat s) \geq  \gamma \frac{v (\infty )}{2} \int _{\widehat s}^{\infty }\int
_s^{\infty} p(t) dt\,ds = \gamma \frac{v (\infty )}{2}\bigg( \int
_{\widehat s}^{\infty} tp(t) dt - \widehat s\int _{\widehat
s}^{\infty} p(t) dt\bigg)\,,$$ thus contradicting \eqref{n112,'}.
\qed
\par\noindent
\bigskip
\par\noindent
{\bf Proof of Theorem \ref{boundlambda}} Assumption \eqref{eiinftylambdabis} implies \eqref{eiinftylambda}, via Fubini's theorem. Hence, the conclusion follows via Theorem \ref{eigencorinf}.
\qed
\smallskip
\par\noindent
{\bf Proof of Theorem \ref{eigensharplambda}} By Corollary
\ref{manifoldsofrevolution}, the same argument as in the proof of
Theorem \ref{eigensharp} provides a manifold $M$, fulfilling
\eqref{1001bis}, but not \eqref{eiinftylambdabis}, on which the
Laplacian has an unbounded eigenfunction. \qed

\section{Applications}\label{appl}
We conclude with applications of our results to two special
instances. Both of them involve families of noncompact manifolds.
However, the former is less pathological, and can be handled
either by isoperimetric or by isocapacitary methods, with the same
output.
That isocapacitary
inequalities can actually yield sharper conclusions than those
obtained via isoperimetric inequalities is
 demonstrated by the latter example, which deals with a class of
manifolds with
a more complicated geometry.
%


\subsection{A family of manifolds of revolution with borderline decay}\label{applrev}
Consider a one-parameter family  of manifolds of revolution $M$
 as in Section \ref{rev}, whose profile $\varphi : [0,
\infty ) \to [0, \infty )$ is such that
\begin{equation}\label{ex1}
\varphi (r) = e^{-r^\alpha } \qquad\quad \hbox{for large $r$,}
\end{equation}
and fulfills the assumptions of Theorem \ref{disisopA} (Figure 1).
\begin{figure}[ht]
\begin{center}

\includegraphics[height=4cm]{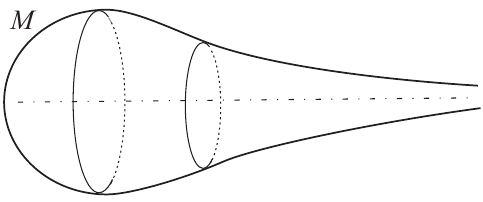}
\end{center}
        \label{figrevolution}
        \caption{A manifold of revolution}
\end{figure}
This theorem  tells us that
\begin{equation}\label{ex2}
\lambda _M(s) \approx s \big(\log (1/s))\big)^{1- 1/\alpha} \qquad
\quad \hbox{near $0$,}
\end{equation}
and
\begin{equation}\label{ex3}
\nu _M(s) \approx \bigg(\int _{s}^{\m2 } \frac{dr}{\lambda _M
(r)^2}\bigg)^{-1} \approx s \big(\log (1/s)\big)^{2- 2/\alpha}
\qquad \hbox{near $0$.}
\end{equation}
An application of Theorem \ref{eigencor} ensures, via \eqref{ex3},
that all eigenfunctions of the Laplacian on $M$ belong to
$L^q(M)$, provided that
\begin{equation}\label{alpha1}
 \alpha >1.
 \end{equation}
   On the other
hand, from Theorem \ref{eigencorinf} and equation \eqref{ex3} one
can infer that the relevant eigenfunctions are  bounded under the
more stringent assumption that
\begin{equation}\label{alpha2}
\alpha
>2.
\end{equation}
%
\par\noindent
The same conclusions can be derived via Theorems
 \ref{Lqlambdabis} and \ref{boundlambda},
respectively.
Thus, as for any other manifold of revolution  of the kind
considered in Theorem \ref{disisopA} (see Corollary
\ref{manifoldsofrevolution}),
 isoperimetric and isocapacitary methods lead to equivalent
results for this family of noncompact manifolds.
\par
Note that, if $\alpha >1$, then, by Proposition
\ref{eigencorconstant}, there exist constants $C_1=C_1(q)$ and
$C_2=C_2(q)$  such that
$$\|u\|_{L^q(M)} \leq C_1 e^{C_2 \gamma ^{\frac{\alpha}{2\alpha
-2}}} \|u\|_{L^2(M)}$$ for any eigenfunction $u$ of the Laplacian
associated with the eigenvalue $\gamma$. Moreover, if $\alpha >2$,
then by Proposition \ref{eigencorinfconstant},
$$\|u\|_{L^\infty(M)} \leq C_1 e^{C_2 \gamma ^{\frac{\alpha}{\alpha
-2}}} \|u\|_{L^2(M)}$$ for some absolute constants $C_1$ and $C_2$
and for every for any eigenfunction $u$ associated with $\gamma$.
In both cases, the existence of eigenfunctions of the Laplacian is
guaranteed by condition \eqref{1001} -- see the comments following
Theorem \ref{eigencor}.


\bigskip
\bigskip
\bigskip

\subsection{A family of manifolds with clustering submanifolds}\label{applch}

Here, we are concerned with a class of non compact surfaces $M$ in
$\R ^3$,  which are patterned on an example appearing in \cite{CH}
and dealing with a planar domain. Their main feature is that they
contain a sequence of mushroom-shaped submanifolds $\{N^k\}$
clustering at some point (Figure 2).
\begin{figure}[ht]
\begin{center}
\includegraphics[height=8cm]{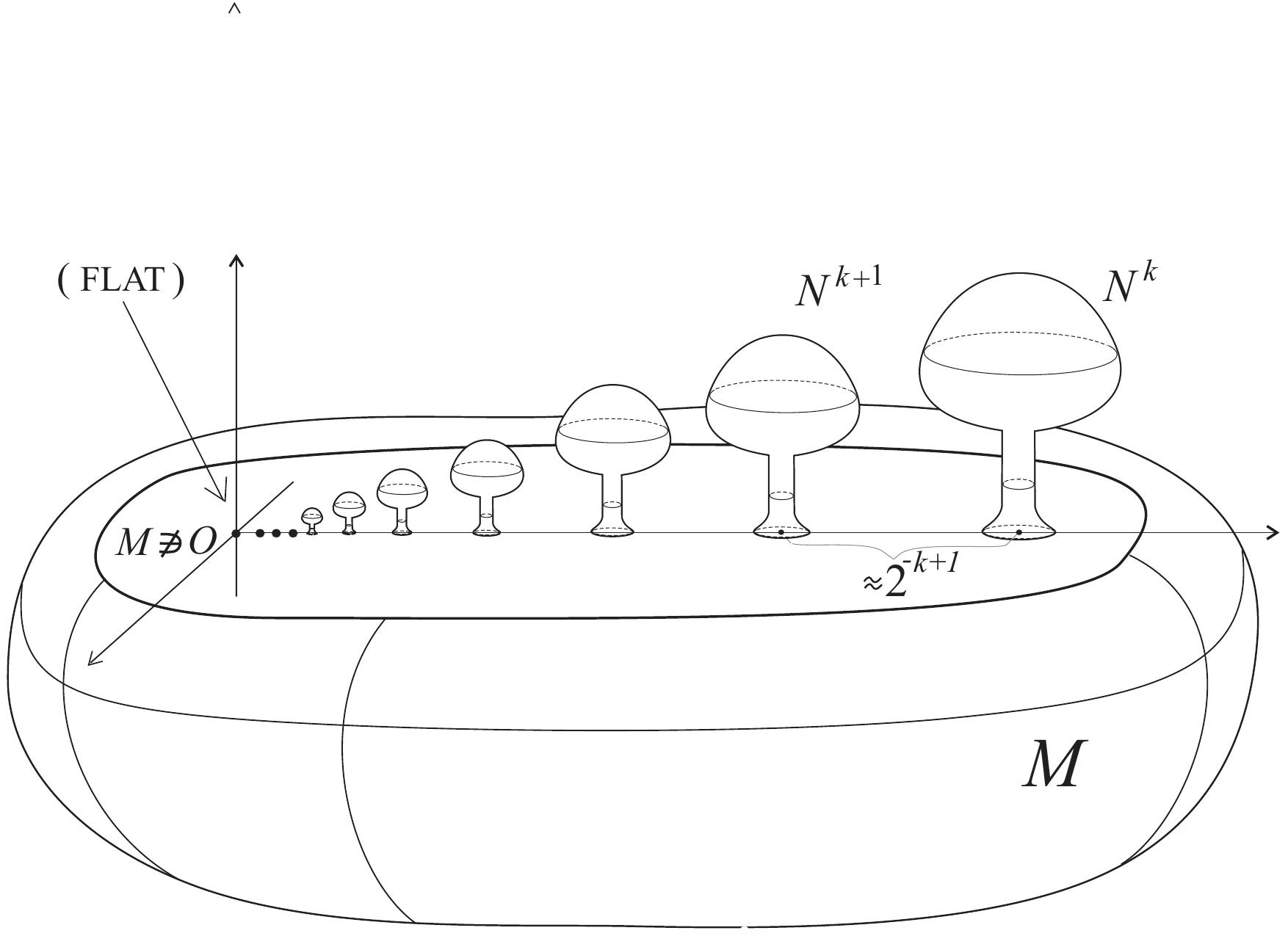}
\end{center}
        \label{figcourant}
        \caption{A manifold
with a family of clustering  submanifolds}
\end{figure}
Let us emphasize that the  submanifolds $\{N^k\}$ are not obtained
just by dilation of each other. Roughly speaking, the diameter of
the head and the length of the neck  of $N^k$ decay to $0$ as
$2^{-k}$ when $k \to \infty$, whereas the width of the neck of
$N^k$ decays to $0$ as $\sigma (2^{-k})$, where $\sigma$ is a
function such that
$$\lim _{s\to 0} \frac{\sigma (s) }{s}=0.$$
 The
isoperimetric and isocapacitary functions of $M$ depend on the
behavior of $\sigma$ at $0$ in a way described in the next result
(Proposition \ref{chcorollary}). Qualitatively, a faster decay to
$0$ of the function $\sigma (s)$ as $s \to 0$ results in a faster
decay to $0$ of $\lambda _M(s)$ and $\nu _M(s)$, and hence in a
 manifold $M$ with a  more irregular geometry. Proposition \ref{chcorollary} is a
special case of Proposition \ref{couranthilbert}, dealing with the
isocapacitary function $\nu _{M,p}$ of $M$ for arbitrary $p \in
[1,2]$, stated and proved below. We also refer to the proof of
Proposition \ref{couranthilbert}  for a more precise definition of
the manifold $M$.

\begin{proposition}\label{chcorollary}
Let $M$ be the $2$-dimensional manifold in Figure 2. Assume that
$\sigma : [0, \infty ) \to [0, \infty )$ is an increasing function
of class $\Delta _2$ such that
\begin{equation}\label{ch-20}
\frac{s^{\beta +1}}{\sigma (s)} \quad \hbox{is non-increasing}
\end{equation}
for some $\beta >0$.
\par\noindent (i) If
\begin{equation*}
\frac{s^{2}}{\sigma (s)} \quad \hbox{is non-decreasing},
\end{equation*}
 then
\begin{equation}\label{ch2}
\lambda _{M} (s) \approx \sigma (s^{1/2})\qquad \hbox{near $0$.}
\end{equation}
(ii) If
\begin{equation*}
\frac{s^{3}}{\sigma (s)} \quad \hbox{is non-decreasing},
\end{equation*}
 then
\begin{equation}\label{ch1}
\nu _{M} (s) \approx \sigma (s^{1/2})s^{-\frac{1}{2}}\qquad
\hbox{near $0$.}
\end{equation}
\end{proposition}

\par
Owing to equation \eqref{ch1}, one can derive the following
conclusions from Theorems  \ref{eigencor} and \ref{eigencorinf},
involving the isocapacitary function $\nu _M$. Assume that
\begin{equation}\label{ex5}
\lim _{s \to 0} \frac{s^3}{\sigma (s)} =0\,.
\end{equation}
Then any eigenfunction of the Laplacian on $M$ belongs to $L^q(M)$
for any $q<\infty$. If \eqref{ex5} is strengthened to
\begin{equation}\label{ex6}
\int _0 \frac{s^2}{\sigma (s)}\,ds < \infty\,,
\end{equation}
then  any eigenfunction of the Laplacian on $M$ is in fact
bounded.
\par\noindent
Conditions \eqref{ex5} and \eqref{ex6} are weaker than parallel
conditions which are obtained from an application of Theorems
\ref{Lqlambdabis} and \ref{boundlambda} and \eqref{ch2},
and read
\begin{equation}\label{ex7}
\lim _{s \to 0} \frac{s^2}{\sigma (s)} =0\,,
\end{equation}
and
\begin{equation}\label{ex8}
\int _0 \frac{s^3}{\sigma (s)^2}\,ds < \infty\,,
\end{equation}
respectively. For instance, if   $b>1$ and $$\sigma (s) = s^b
\quad \quad \hbox{for $s >0$,}$$
  then \eqref{ex5} and \eqref{ex6}
amount to  $b <3$, whereas \eqref{ex7} and \eqref{ex8} are
equivalent to the more stringent condition that  $b <2$.
\par\noindent
Since, by \eqref{ch1}, $\nu _M(s) \approx s^{\frac{b-1}2}$, from
Examples \ref{example1} and \ref{example2}
%
we deduce that there exists a constant
$C=C(q)$ such that
$$\|u\|_{L^q(M)} \leq C \gamma ^{\frac{q-2}{q(3-b)} }\|u\|_{L^2(M)}$$
for every $q \in (2, \infty ]$ and for any eigenfunction $u$ of
the Laplacian associated with the eigenvalue $\gamma$. The
existence of such eigenfunction follows from condition
\eqref{1001}, as explained in the comments following Theorem
\ref{eigencor}.

%
%
%

\begin{proposition}\label{couranthilbert}
Let $M$ be the $2$-dimensional manifold in Figure 2. Let $1 \leq p
\leq 2$, and let $\sigma : [0, \infty ) \to [0, \infty )$ be an
increasing function of class $\Delta _2$. Then there exist a
constant $C$ such that
\begin{equation}\label{ch-3}
\nu _{M,p} (s) \leq C \sigma (s^{1/2})s^{-\frac{p-1}{2}} \qquad
\hbox{near $0$.}
\end{equation}
Assume, in addition, that
\begin{equation}\label{ch-2}
\frac{s^{\beta +1}}{\sigma (s)} \quad \hbox{is non-increasing}
\end{equation}
for some $\beta >0$, and
\begin{equation}\label{ch-1}
\frac{s^{p+1}}{\sigma (s)} \quad \hbox{is non-decreasing}.
\end{equation}
Then
\begin{equation}\label{ch3}
\nu _{M,p} (s) \approx \sigma (s^{1/2})s^{-\frac{p-1}{2}} \qquad
\hbox{near $0$.}
\end{equation}
\end{proposition}

Note that equation \eqref{ch2} of Proposition \ref{chcorollary}
follows from \eqref{ch3}, owing to property \eqref{nu1lambda},
whereas equation \eqref{ch1} agrees with \eqref{ch3} for $p=2$.
\par
 One step in the proof of Proposition \ref{couranthilbert} makes use of Orlicz spaces. Recall that given a non-atomic,
$\sigma$-finite measure space $(\mathcal R , m)$ and a Young
function $A$, namely a convex function
 from $ [0,
\infty )$ into $[0, \infty )$ vanishing at $0$, the Orlicz space
$L^A(\mathcal R)$ is  the Banach space of those real-valued
$m$-measurable functions $f$ on $\mathcal R$ whose Luxemburg norm
%
%
\begin{equation*}
 \|f\|_{L^A(\mathcal R)}= \inf \left\{ \lambda >0 :  \int_{\mathcal R}A
\left( \frac{|f|}{\lambda} \right) dm \leq 1 \right\}\,,
\end{equation*}
is finite. A generalized H\"older type inequality in Orlicz spaces
tells us that if $A_i$, $i=1,2,3$, are Young functions such that
$A_1^{-1}(r)A_2^{-1}(r)\leq C A_3^{-1}(r)$, then there exists a
constant $C'$ such that
\begin{equation}\label{holder}
\|fg\|_{L^{A_3}(\mathcal R)} \leq C' \|f\|_{L^{A_1}(\mathcal R)}
\|g\|_{L^{A_2}(\mathcal R)}
\end{equation}
for every $f \in L^{A_1}(\mathcal R)$ and $g \in L^{A_2}(\mathcal
R)$ \cite{Oneil}.

\medskip
\par\noindent
{\bf Proof of Proposition \ref{couranthilbert}}.
%


\par\noindent
{\bf Part I}.
Here we show that, if  \eqref{ch-2} and
 \eqref{ch-1} are in force, then there exists a constant $C$ such
that
\begin{equation}\label{ch4}
\nu _{M,p} (s) \geq C \sigma (s^{1/2})s^{-\frac{p-1}{2}} \qquad
\hbox{for $s \in (0, \m2 )$.}
\end{equation}
We split the proof of \eqref{ch4} is steps.
\par\noindent
{\bf Step 1}. Fixed $\varepsilon _0 >0$ and $\varepsilon \in (0,
\varepsilon _0)$, let $N_\varepsilon = Q \cup P_\varepsilon \cup
R_\varepsilon$ be the auxiliary surface of revolution in $\R^3$
depicted in Figure 3.

%
%

\begin{figure}[ht]
\begin{center}
\includegraphics[height=8cm]{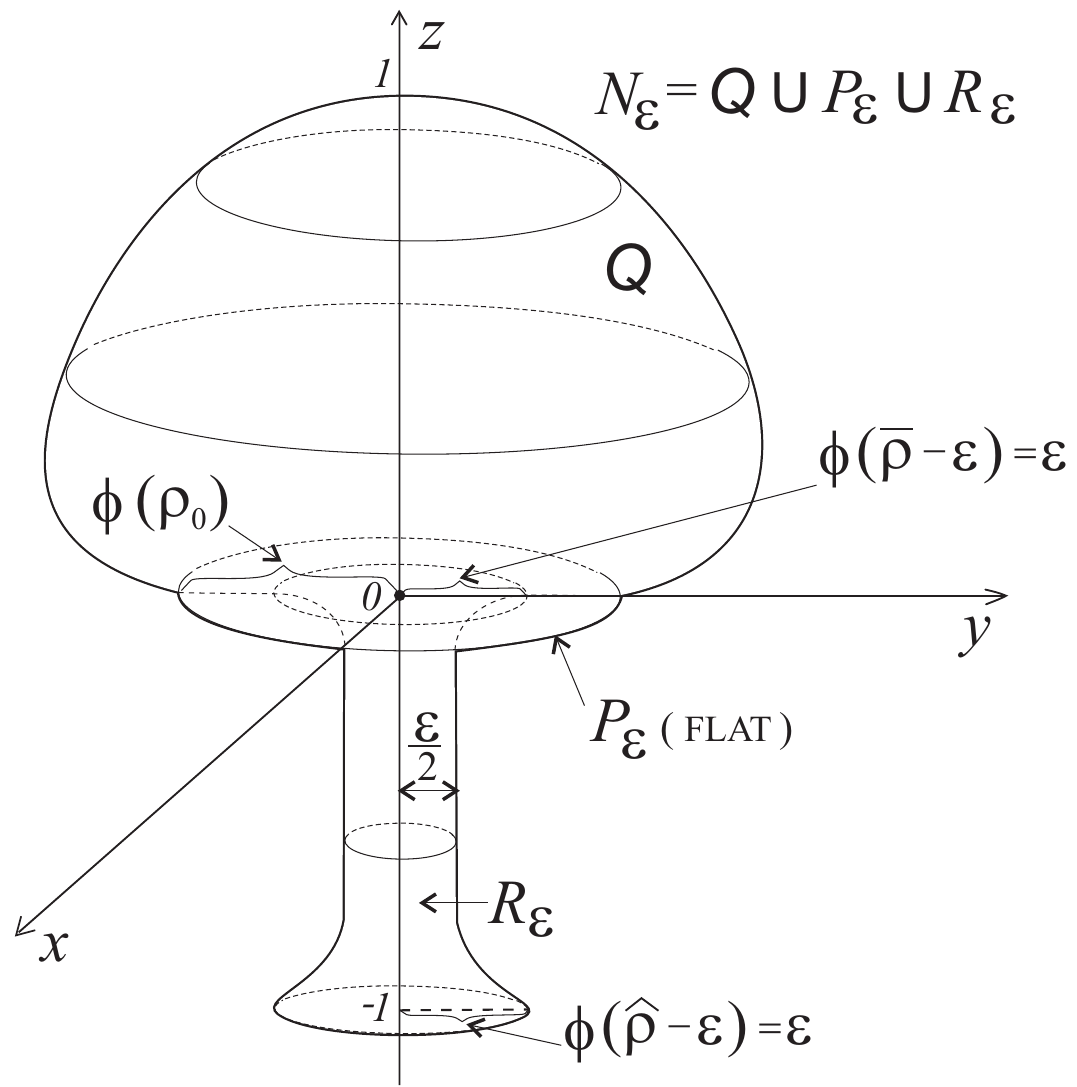}
\end{center}
        \label{figauxiliary}
        \caption{An auxiliary
submanifold}
\end{figure}

Let $q= \frac {2p}{2-p}$ if $p<2$, and let $q$ be a sufficiently
large number, to be chosen later, if $p=2$. We shall show that
\begin{equation}\label{ch4'}
\bigg(\int _{Q \cup P_\varepsilon }|u|^{q}d\hdue \bigg)^{\frac
p{q}} \leq \frac{C}{\varepsilon}\bigg(\int _{N_\varepsilon}|\nabla
u|^p d\hdue + \int _{\partial N_\varepsilon} |u|^p d\huno \bigg)
\end{equation}
for every $u \in  W^{1,p}(N_\varepsilon )$, and for some  constant
$C$ independent of $\varepsilon$ and $u$.
\par\noindent
Let $(\rho, \vartheta ) \in [0, \hat{\rho} - \varepsilon ) \times
[0, 2\pi ]$ be  geodesic coordinates on $N_\varepsilon$ with
respect to the point $(0,0,1)$, and
\begin{equation*}
\begin{cases}
x & =  \phi (\rho )\cos \,\vartheta \\
y & = \phi (\rho )\sin \, \vartheta  \\
z & = \psi (\rho )
\end{cases}
\end{equation*}
%
be a parametrization of $N_\varepsilon $ for some  given smooth
functions $\phi , \psi : [0, \hat \rho - \varepsilon ] \to [0,
\infty )$. In particular, $ \phi '(\rho )^2 + \psi '(\rho )^2=1$
for $\rho \in [0, \hat{\rho} - \varepsilon )$.
 The functions $\phi$ and $\psi$ are independent of
$\varepsilon$ in $[0, \rho _0]$ and (up to a translation, of
lenght $\varepsilon$, in the variable $\rho$ ) in
$[\overline{\rho} - \varepsilon, \hat{\rho} - \varepsilon ]$; on
the other hand, since $P_\varepsilon$ is a flat annulus, we have
that $\psi (\rho) =0$ for $\rho \in [\rho _0, \overline{\rho}
-\varepsilon ]$ and that $\phi $ is an affine function in the same
interval.
\par\noindent Thus, the metric on $M$ is given by
$$d s^2 = d \rho ^2 + \phi (\rho )^2 d \vartheta ^2\,.$$
%
%
  In particular,
$$ \int _{N_\varepsilon} f d \hdue = \int _0^{2\pi } \int _0^{\hat{\rho} - \varepsilon
} f \, \phi (\rho )\, d\rho d\vartheta $$ for any integrable
function $f$ on $M$.
Moreover, if $u \in W^{1,p}(N_\varepsilon )$,
\begin{equation}\label{ch8}
|\nabla u| = \sqrt{u_\rho^2 + \frac {u_\vartheta ^2}{\phi (\rho)^2
}} \qquad \hbox{a.e. in $N_\varepsilon$}\,.
\end{equation}
Define
\begin{equation}\label{ch9}
\overline{u}(\vartheta ) = \frac 1{\int _0^{\overline{\rho}
-\varepsilon}\phi(\rho )d\rho} \int _0^{\overline{\rho}
-\varepsilon} u(\rho , \vartheta ) \phi(\rho )d\rho  \qquad \hbox{for a.e.  $\vartheta \in [0, 2\pi ]$}\,.
\end{equation}
One has that
\begin{equation}\label{ch10}
\bigg(\int _{Q \cup P_\varepsilon }|u|^{q}d\hdue \bigg)^{\frac
p{q}} \leq 2^{p-1} \bigg[\bigg(\int _{Q \cup P_\varepsilon
}|u-\overline{u}|^{q}d\hdue \bigg)^{\frac p{q}}+ \bigg(\int _{Q
\cup P_\varepsilon }|\overline{u}|^{q}d\hdue \bigg)^{\frac
p{q}}\bigg]
\,,
\end{equation}
where $\overline u$ is regarded as a function defined on
$N_\varepsilon$. It is easily verified that the function $u-
\overline u$ has mean value $0$ on $Q \cup P_\varepsilon$. Thus,
by a standard Poincar\'e inequality,
\begin{equation}\label{ch11}
\bigg(\int _{Q \cup P_\varepsilon
}|u-\overline{u}|^{q}d\hdue \bigg)^{\frac p{q}}\leq C \int _{Q \cup P_\varepsilon }|\nabla u|^p d\hdue\,,
\end{equation}
for some constant $C$ independent of $\varepsilon$ and $u$.
This is a
consequence of the fact that $Q$ is independent of $\varepsilon$,
and $P_\varepsilon $ is an open subset of $\R ^2$ (an annulus)
enjoying the cone property for some cone independent of
$\varepsilon$.
\par\noindent
Next, the following chain holds:
\begin{align}\label{ch12}
\bigg(\int _{Q \cup  P_\varepsilon }& |\overline{u}|^{q}d\hdue
\bigg)^{\frac p{q}}  = \bigg(\int _0^{2\pi}\int _0^{\overline
\varrho - \varepsilon}\bigg| \frac 1{\int _0^{\overline{\rho}
-\varepsilon} \phi (r )dr} \int _0^{\overline{\rho} -\varepsilon}
u(r , \vartheta ) \phi (r )\, dr
\bigg|^{q} \phi (\rho )d\rho d\vartheta\bigg)^{\frac p{q}}\\
\nonumber & = \bigg( \int _0^{\overline{\rho} -\varepsilon} \phi
(\rho )d\rho\bigg)^{\frac p{q}} \bigg(\int _0^{2\pi}\bigg| \frac
1{\int _0^{\overline{\rho} -\varepsilon} \phi (r )dr} \int
_0^{\overline{\rho} -\varepsilon} u(r , \vartheta ) \phi (r )dr
\bigg|^{q} d\vartheta\bigg)^{\frac p{q}} \\ \nonumber & \leq
\bigg( \int _0^{\overline{\rho} -\varepsilon} \phi (\rho
)d\rho\bigg)^{\frac p{q}}(2\pi)^{\frac p{q}}\sup _{\vartheta \in
[0, 2\pi ]}\bigg| \frac 1{\int _0^{\overline{\rho} -\varepsilon}
\phi (r )dr} \int _0^{\overline{\rho} -\varepsilon} u(r ,
\vartheta ) \phi (r )dr \bigg|^{p}
\\ \nonumber & \leq C \bigg( \int
_0^{\overline{\rho} -\varepsilon} \phi (\rho )d\rho\bigg)^{\frac
p{q}}\bigg[\bigg(\int _0^{2\pi}\frac 1{\int _0^{\overline{\rho}
-\varepsilon} \phi (r )dr} \int _0^{\overline{\rho} -\varepsilon}
|u_\vartheta(r , \vartheta )| \phi (r )dr  d\vartheta\bigg)^p
\\ \nonumber & \qquad + \bigg(\frac 1{2\pi} \int _0^{2\pi}\frac 1{\int _0^{\overline{\rho}
-\varepsilon} \phi (r )dr} \int _0^{\overline{\rho} -\varepsilon}
|u(r , \vartheta )|\phi (r )dr d\vartheta\bigg)^p\bigg]
\\ \nonumber &
\leq C \bigg( \int _0^{\overline{\rho} -\varepsilon} \phi (\rho
)d\rho\bigg)^{\frac p{q}}\bigg[\bigg(\int _0^{2\pi}\frac 1{\int
_0^{\overline{\rho} -\varepsilon} \phi (r )dr} \int
_0^{\overline{\rho} -\varepsilon} |\nabla u| \phi (r)^2 dr
d\vartheta\bigg)^p
\\ \nonumber & \qquad + \bigg(\frac 1{2\pi} \int _0^{2\pi}\frac 1{\int _0^{\overline{\rho}
-\varepsilon} \phi (r )dr} \int _0^{\overline{\rho} -\varepsilon}
|u(r , \vartheta )|\phi (r )dr d\vartheta\bigg)^p\bigg]
\\ \nonumber &
\leq  C \bigg( \int _0^{\overline{\rho} -\varepsilon} \phi (\rho
)d\rho\bigg)^{\frac p{q}}\bigg[\max |\phi|^p \bigg(\int
_0^{2\pi}\frac 1{\int _0^{\overline{\rho} -\varepsilon} \phi (r
)dr} \int _0^{\overline{\rho} -\varepsilon} |\nabla u| \phi (r )dr
d\vartheta\bigg)^p
 \\ \nonumber & \qquad + \bigg(\frac 1{2\pi} \int
_0^{2\pi}\frac 1{\int _0^{\overline{\rho} -\varepsilon} \phi (r
)dr} \int _0^{\overline{\rho} -\varepsilon} |u(r , \vartheta
)|\phi (r )dr d\vartheta\bigg)^p\bigg]
\end{align}
\begin{align*} &
 \leq C \bigg( \int _0^{\overline{\rho} -\varepsilon} \phi
(\rho )d\rho\bigg)^{\frac p{q}}\bigg[\max |\phi|^p\bigg(\int
_0^{2\pi}\bigg(\frac 1{\int _0^{\overline{\rho} -\varepsilon} \phi
(r )dr} \int _0^{\overline{\rho} -\varepsilon} |\nabla u|^p \phi
(r )dr \bigg)^{1/p} d\vartheta\bigg)^p \\ \nonumber
 & \qquad + \bigg(\frac 1{2\pi} \int _0^{2\pi}\bigg(\frac 1{\int
_0^{\overline{\rho} -\varepsilon} \phi (r )dr} \int
_0^{\overline{\rho} -\varepsilon} |u(r , \vartheta )|^p\phi (r
)dr\bigg)^{1/p} d\vartheta\bigg)^p\bigg]
\\ \nonumber
 &
\leq  C \bigg( \int _0^{\overline{\rho} -\varepsilon} \phi (\rho
)d\rho\bigg)^{\frac p{q}}\bigg[\frac {\max |\phi|^p(2 \pi
)^{p-1}}{\int _0^{\overline{\rho} -\varepsilon} \phi (r )dr} \int
_0^{2\pi}\int _0^{\overline{\rho} -\varepsilon} |\nabla u|^p \phi
(r )dr  d\vartheta \\ \nonumber & \qquad + \frac 1{2\pi \int
_0^{\overline{\rho} -\varepsilon} \phi (r )dr} \int _0^{2\pi}\int
_0^{\overline{\rho} -\varepsilon} |u(r , \vartheta )|^p\phi (r )dr
d\vartheta \bigg]
\\ \nonumber & \leq C' \bigg[\int _{Q\cup P_\varepsilon} |\nabla
u|^p d\hdue + \int _0^{2\pi}\int _0^{\overline{\rho} -\varepsilon}
|u(r , \vartheta )|^p\phi (r )dr d\vartheta\bigg],
\end{align*}
for some constants $C$ and $C'$ independent of $\varepsilon$ and
$u$. Note that a rigorous derivation of the inequality between the
leftmost and rightmost sides of \eqref{ch12} requires an
approximation argument of $u$ by smooth functions. Since, for a.e.
$\vartheta \in [0, 2 \pi ]$,
\begin{equation}\label{ch13}
u(\rho ,\vartheta ) = u(\hat{\rho} - \varepsilon, \vartheta ) -
\int ^{\hat{\rho} - \varepsilon}_\rho u_\rho (r , \vartheta )dr\,
\qquad \quad \hbox{for $\rho \in (0, \hat \rho -\varepsilon )$,}
\end{equation}
we have that
\begin{equation}\label{ch14}
|u(\rho ,\vartheta )|^p \leq C |u(\hat{\rho} - \varepsilon,
\vartheta )|^p+ C \int ^{\hat{\rho} - \varepsilon}_\rho |u_\rho (r
, \vartheta )|^pdr \qquad \hbox{for a.e. $(\rho , \vartheta ) \in
(0 , \widehat \rho - \varepsilon ) \times (0, 2\pi )$\,,}
\end{equation}
for some constant $C$ independent of $\varepsilon$ and $u$. Thus,
\begin{align}\label{ch15}
\int _0^{2\pi}& \int _0^{\overline{\rho} -\varepsilon} |u(\rho ,
\vartheta )|^p\phi (\rho )d\rho d\vartheta \\ \nonumber & \leq C
\int _0^{2\pi}\int _0^{\hat{\rho} -\varepsilon} \bigg(\int
^{\hat{\rho} - \varepsilon}_\rho |u_\rho (r , \vartheta
)|^pdr\bigg) \phi (\rho ) d\rho d\vartheta + C \int _0^{2\pi}\int
_0^{\hat{\rho} -\varepsilon} |u(\hat{\rho} - \varepsilon,
\vartheta )|^p \phi (\rho )d\rho d\vartheta
\\ \nonumber & \leq
C \int _0^{2\pi}\int _0^{\hat{\rho} -\varepsilon} \bigg(\int
^{\hat{\rho} - \varepsilon}_\rho |\nabla u (r, \vartheta)|^p
dr\bigg) \phi (\rho ) d\rho d\vartheta + C \int _0^{2\pi}\int
_0^{\hat{\rho} -\varepsilon} |u(\hat{\rho} - \varepsilon,
\vartheta )|^p \phi (\rho )d\rho d\vartheta
\\ \nonumber & =
C \int _0^{2\pi}\int _0^{\hat{\rho} -\varepsilon}|\nabla u(r,
\vartheta) |^p \bigg(\int ^{r}_0 \phi (\rho )  d\rho \bigg) dr
d\vartheta + C \bigg(\int _0^{\hat{\rho} -\varepsilon}\phi (\rho
)d\rho\bigg)\bigg(\int _0^{2\pi} |u(\hat{\rho} - \varepsilon,
\vartheta )|^p d\vartheta\bigg)
\\ \nonumber & \leq C \bigg(\sup _{r \in (0, \hat{\rho}
-\varepsilon)}\frac{\int ^{r}_0 \phi (\rho ) d\rho }{\phi (r
)}\bigg) \int _0^{2\pi}\int _0^{\hat{\rho} -\varepsilon}|\nabla
u(r, \vartheta) |^p  \phi (r ) dr   d\vartheta  \\
\nonumber & \quad + C \bigg(\int _0^{\hat{\rho} -\varepsilon}\phi
(\rho )d\rho\bigg)\bigg(\int _0^{2\pi} |u(\hat{\rho} -
\varepsilon, \vartheta )|^p d\vartheta \bigg)
\\ \nonumber & \leq C' \frac{1}{\phi ((\hat{\rho} - \overline{\rho})/2))}
 \int _0^{2\pi}\int _0^{\hat{\rho}
-\varepsilon}|\nabla u(r, \vartheta) |^p  \phi (r ) dr d\vartheta
+ C \bigg(\int _0^{\hat{\rho} -\varepsilon}\phi (\rho )d\rho\bigg)
\bigg(\int _0^{2\pi} |u(\hat{\rho} - \varepsilon, \vartheta )|^p
d\vartheta\bigg)
\\ \nonumber & = C' \frac{1}{\varepsilon /2}
 \int _{N_\varepsilon}|\nabla u |^p  d\hdue +  \frac{C''}{\varepsilon}\int
_{\partial N_\varepsilon} |u|^p d\huno
\,,
\end{align}
for some constants $C$, $C'$ and $C''$ independent of
$\varepsilon$ and $u$. Note that the last inequality holds since
$\phi$ is increasing in a neighborhood of $0$. \par\noindent
Inequality \eqref{ch4'} follows from \eqref{ch10}, \eqref{ch11},
\eqref{ch12} and \eqref{ch15}.
\par\noindent
{\bf Step 2}. Let $N_{\varepsilon , \delta}$ be the manifold
obtained on scaling $N_\varepsilon$ by a factor $\delta$. Thus,
$N_{\varepsilon , \delta}$ is parametrized  by $(x,y,z) = (
\delta\phi (\rho  )\cos \vartheta , \delta\phi (\rho  )\sin
\vartheta , \delta\psi (\rho ))$.
%
%
%
%
\par\noindent
Let $u \in W^{1,p}(N_{\varepsilon , \delta})$. From
\eqref{ch4'} we obtain that
\begin{equation}\label{ch16}
\delta ^{-\frac {2p}{q}}\bigg(\int _{Q_\delta \cup P_{\varepsilon
, \delta}}|u|^{q}d\hdue \bigg)^{\frac p{q}} \leq \frac{C}{
\varepsilon}\bigg(\delta ^{p-2}\int _{N_{\varepsilon ,
\delta}}|\nabla u|^p d\hdue  + \delta ^{-1}\int _{\partial
N_{\varepsilon , \delta}} |u|^p d\huno \bigg),
\end{equation}
where $Q_\delta$ and $P_{\varepsilon , \delta}$ denote the subsets
of $N_{\varepsilon , \delta}$ obtained on scaling $Q$ and
$P_{\varepsilon }$, respectively, by a factor $\delta$.
Now, let $A$ be a Young function whose inverse
satisfies
\begin{equation}\label{ch17}
A^{-1}(\delta ^{-2}) \approx \frac{\delta ^{p-1}}{\sigma (\delta)}
\qquad \quad \hbox{for $\delta >0$}.
\end{equation}
Notice that such a function $A$ does exist. Indeed,  the function
$H : (0, \infty ) \to [0, \infty )$ given by $H(t)=
\frac{t^{-\frac{p-1}{2}}}{\sigma (t^{-\frac 12})}$ for $t >0$ is
increasing  by \eqref{ch-2}, and the function $\frac {H(t)}{t}$ is
non-increasing by \eqref{ch-1}. Thus, $\frac{H^{-1}(\tau )}{\tau}$
is a non-decreasing function, and, on taking
$$A(t) = \int _0^t \frac{H^{-1}(\tau )}{\tau } \, d \tau \qquad \hbox{for $t \geq 0$,}$$
equation \eqref{ch17} holds, inasmuch as $A(t) \approx H^{-1}(t)$
for $t \geq 0$.
%
%
%
Next, we claim that  a Young function $E$ exists whose inverse
fulfils
\begin{equation}\label{ch18}
E^{-1}(\tau) \approx \frac{A^{-1}(\tau)}{\tau ^{p/q}} \qquad
\hbox{for $\tau
> 0$.}
\end{equation}
To see this, note that  the function
$J(\tau)=\frac{A^{-1}(\tau)}{\tau ^{p/q}}$ is equivalent to an
increasing function $F(\tau)$ (for sufficiently large $q$,
depending on $\beta$, if $p=2$) by \eqref{ch-2}, and that the
function $\frac{J(\tau)}{\tau} = \frac{A^{-1}(\tau)}{\tau
^{1+p/q}}$ is trivially decreasing. Set $J_1(\tau ) = \frac{J(\tau
)}{\tau}$. Thus, $\frac {F(\tau )}\tau \approx J_1(\tau)$ for
$\tau >0$. As a consequence, one can show that $\frac {F^{-1}(t
)}{t } \approx \frac 1{J_1 (F^{-1}(t ))}$, an increasing function.
Thus the function $E$ given by
$$E(t ) = \int _0^t \frac {d\tau }{J_1 (F^{-1}(\tau ))}\qquad \hbox{for
$t \geq 0$,}$$ is a Young function, and since $E(t ) \approx
\frac{t}{J_1 (F^{-1}(t ))} \approx F^{-1}(t)$, one has that
$E^{-1}(\tau) \approx F(\tau) \approx J(\tau)
=\frac{A^{-1}(\tau)}{t^{p/q}}$, whence \eqref{ch18} follows.
\par\noindent
Owing to \eqref{ch18}, inequality \eqref{holder} ensures that
\begin{multline}\label{ch19}
\||u|^p\|_{L^A(Q_\delta \cup P_{\varepsilon , \delta})} \leq C
\||u|^p\|_{L^{q/p}(Q_\delta \cup P_{\varepsilon ,
\delta})}\|1\|_{L^E(Q_\delta \cup P_{\varepsilon , \delta})} \\= C
\|u\|_{L^{q}(Q_\delta \cup P_{\varepsilon , \delta})}^p\frac
{1}{E^{-1}(1/\hdue (Q_\delta \cup P_{\varepsilon , \delta}))} \leq
C \|u\|_{L^{q}(Q_\delta \cup P_{\varepsilon , \delta})}^p\frac
{1}{E^{-1}(C'/\delta ^2)}\,,
\end{multline}
for some constants $C$ and $C'$ independent of $\varepsilon$,
$\delta$ and $u$. Combining \eqref{ch16}--\eqref{ch19} yields
\begin{align}\label{ch20}
\||u|^p\|_{L^A(Q_\delta \cup P_{\varepsilon , \delta})} & \leq
\frac {C\sigma (\delta )}{\varepsilon \delta } \int
_{N_{\varepsilon , \delta}}|\nabla u|^p d\hdue + \frac {C\sigma
(\delta )}{\varepsilon\delta ^{p}}\int _{\partial N_{\varepsilon ,
\delta}} |u|^p d\huno
 \,.
\end{align}
Now, choose
$$ \varepsilon = \frac {\sigma (\delta )}{\delta}\,,$$
and obtain from \eqref{ch20}
\begin{equation}\label{ch21}
\||u|^p\|_{L^A(Q_\delta \cup P_{\frac {\sigma (\delta )}{\delta}
,\delta })} \leq C\int _{N_{\frac {\sigma (\delta )}{\delta},
\delta}}|\nabla u|^p d\hdue + C\delta ^{1-p}\int _{\partial
N_{\frac {\sigma (\delta )}{\delta}, \delta}} |u|^p d\huno\,.
\end{equation}
Note also that
\begin{equation}\label{riscal}
\huno (\partial N_{\frac {\sigma (\delta )}{\delta}, \delta}) =
2\pi \sigma (\delta ),
\end{equation}
a fact that will be tacitly used in what follows.
\par\noindent
{\bf Step 3}. Choose $\delta _k = 2^{-k}$ for $k \in \N$, and
denote
$$Q^k =Q_{\delta _k}, \quad P^k =P_{\frac {\sigma (\delta _k
)}{\delta _k} ,\delta _k}, \quad  N^k= N_{\frac {\sigma (\delta _k
)}{\delta _k}, \delta _k}.$$
 Define the manifold $M$ in such a way
that the distance between the centers of the circumferences
$\partial N^k$  and $\partial N^{k+1}$ equals $2^{-k+1}$. Given $u
\in W^{1,p}(M)$, one has that
\begin{equation}\label{ch22}
\||u|^p\|_{L^A(\cup _k (Q^k\cup P^k)} \leq \sum _{k}
\||u|^p\|_{L^A( Q^k\cup P^k)}\,,
\end{equation}
and
\begin{equation}\label{ch23}
\int _{\cup _k N^k}|\nabla u|^p d\hdue = \sum _{k \in \N} \int _{
N^k}|\nabla u|^p d\hdue\,.
\end{equation}
Now, notice that the manifold $M$ is flat in a neighborhood of
$\cup _k N^k$.  For $k \in \mathbb \N$, let us denote by $\Omega
_k$ the open set on $M$ bounded by the circumference $\partial
N^k$ (having radius $\sigma (\delta _k)$) and by the boundary of
the square on $M$, with sides parallel to the coordinate axes,
whose side-length is $3\sigma (\delta _k)$, and whose center
agrees with the center of  $\partial N^k$.
Hence, in particular,
\begin{equation}\label{I}
\hdue (\Omega _k) \leq 9 \sigma (\delta _k) ^2.
\end{equation}
Observe that
\begin{equation}\label{A}
\int _{\partial N^k}|u|^p\, d\huno \leq C \sigma (\delta
_k)^{p-1}\int _{\Omega _k}|\nabla u|^pd \hdue + C \sigma (\delta
_k)^{-1}\int _{\Omega _k}|u|^p d \hdue\,
\end{equation}
for some constant $C$ independent of $k$ and $u$. Inequality
\eqref{A} can be derived via a scaling argument applied to a
standard trace inequality for subsets of $\rn$ with a Lipschitz
boundary. Thus,
 \begin{align}\label{B'}
\sum _{k \in \N} \delta _k ^{1-p}& \int _{\partial N^k}|u|^p\,
d\huno \\ \nonumber & \leq C \sum _{k \in \N} \delta _k ^{1-p}
\sigma (\delta _k)^{p-1}\int _{\Omega _k}|\nabla u|^pd \hdue + C
\sum _{k \in \N} \delta _k ^{1-p}  \sigma (\delta _k)^{-1}\int
_{\Omega _k}|u|^p d \hdue\,.
\end{align}
Assumption \eqref{ch-2} ensures that
\begin{equation*}
\lim _{\delta \to 0} \frac {\sigma (\delta )}{\delta } =0\,,
\end{equation*}
and hence, in particular, there exists a constant $C$ such that
\begin{equation}\label{D}
\frac {\sigma (\delta )}{\delta } \leq C \qquad \quad \hbox{if $0
<\delta \leq 1$.}
%
\end{equation}
Consequently,
\begin{equation}\label{B}
 \sum _{k \in \N} \delta _k ^{1-p}  \sigma (\delta _k)^{p-1}\int _{\Omega _k}|\nabla u|^pd \hdue \leq C \int _M |\nabla u|^pd \hdue
\end{equation}
for some constant $C$. As far as the second addend on the
right-hand side of \eqref{B'} is concerned, if $1 \leq p <2$ by
H\"older's inequality and \eqref{I} one has that
\begin{align}\label{F}
\sum _{k \in \N} \delta _k ^{1-p} & \sigma (\delta _k)^{-1}\int
_{\Omega _k}|u|^p d \hdue = \int _M \sum _{k \in \N} \chi _{\Omega
_k} \delta _k ^{1-p}  \sigma (\delta _k)^{-1}|u|^p d \hdue
\\ \nonumber & \leq \bigg(\int _{\cup _k \Omega _{k}}|u|^{\frac{2p}{2-p}}d \hdue \bigg)^{\frac {2-p}{2}}
\Big(9\sum _{k \in \N} \sigma (\delta _k)^2 \big(\delta _k^{1-p}
\sigma (\delta _k )^{-1}\big)^{\frac 2p}\Big)^{\frac p2}
\\ \nonumber & \leq C \bigg(\int _{\cup _k \Omega _{k}}|u|^{\frac{2p}{2-p}}d \hdue \bigg)^{\frac {2-p}{2}}
\bigg(\int _0^1 \frac {\sigma (\delta )^{2-\frac 2p}}{\delta ^{3-
\frac 2p}} \, d \delta\bigg)^{\frac p2}\,,
\end{align}
for some constant $C$ independent of $k$ and $u$. If $p =2$, then
given $a>1$ one similarly has that
\begin{align}\label{G}
\sum _{k \in \N} \delta _k ^{-1} & \sigma (\delta _k)^{-1}\int
_{\Omega _k}|u|^2 d \hdue
 \leq C \bigg(\int _{\cup _k \Omega _{k}}|u|^{2a}d \hdue \bigg)^{\frac {1}{a}}
\bigg(\int _0^1 \frac {\sigma (\delta )^{2-a'}}{\delta ^{1+a'}} \, d \delta\bigg)^{\frac 1{a'}}\,.
\end{align}
Thanks to \eqref{ch-2}, $\int _0^1 \frac {\sigma (\delta )^{2-\frac 2p}}{\delta ^{3- \frac 2p}} \, d \delta < \infty$ if $1 \leq p <2$, and $\int _0^1 \frac {\sigma (\delta )^{2-a'}}{\delta ^{1+a'}} \, d \delta < \infty$ if $p=2$, provided that $a$ is sufficiently large.
\par\noindent
On the other hand, by our choice of $\delta _k$ and of the
distance between the centers of $\partial N^k$ and $\partial
N^{k+1}$, any regular  neighborhood  of $\cup _k
\partial N^k$ in $M$, containing $\cup _k \Omega _k$,  is
a planar domain  having the cone property. Hence, by the Sobolev
inequality, if $1 \leq p <2$
\begin{align*}
 \bigg(\int _{\cup _k \Omega _{k}}& |u|^{\frac{2p}{2-p}}d \hdue \bigg)^{\frac {2-p}{2}} \leq
C \int _{M} \big(|\nabla u|^p + |u|^p
\big)d\hdue \,,
\end{align*}
and, if $p=2$,
$$\bigg(\int _{\cup _k \Omega _{k}}|u|^{2a}d \hdue \bigg)^{\frac {1}{a}}
\leq C
 \int _{M} \big(|\nabla u|^2 + |u|^2
\big)d\hdue
$$
for some constant $C$ independent of $u$. Altogether, we infer
that there exists a constant $C$ such that

\begin{equation}\label{ch31'}
\sum _{k \in \N} \delta _k^{1-p} \sigma (\delta _k)^{-1}\int
_{\partial N^k} |u|^p d\huno \leq  C\bigg( \int _{M} |\nabla u|^p
dx + \int _{M}|u|^{p} dx\bigg)\,.
\end{equation}
Combining \eqref{ch21}, \eqref{ch22}, \eqref{ch23}  and
\eqref{ch31'} tells us that there exists a constant $C$ such that
\begin{equation}\label{ch32}
\||u|^p\|_{L^A(\cup _k (Q^k\cup P^k))}^{1/p} \leq C \big(\|\nabla
u\|_{L^p(M)} + \| u\|_{L^p(M)}\big)
\end{equation}
for every $u \in W^{1,p}(M)$.
\par\noindent
{\bf Step 4}. Denote by
 $R_{\frac {\sigma (\delta )}{\delta},\delta}$ the manifold obtained  on scaling $R_{\frac {\sigma
(\delta )}{\delta }}$ by the factor $\delta$. We shall show that
inequality \eqref{ch32} continues to hold if $\cup _k (Q^k\cup
P^k)$ is replaced by $\cup _k R^k$, where
$$R^k = R_{\frac {\sigma (\delta _k)}{\delta _k},
\delta _k}.$$
 Let $\rho_i$,
$i=1, \dots , m$, be such that $\rho _1= \overline \rho -
\varepsilon$, $\rho _m= \hat \rho - \varepsilon$, the difference
$\rho _{i+1} - \rho _i$ is independent of $i$ for  $i=1, \dots ,
m-1$, and
\begin{equation}\label{approx10}
1 \leq \frac{\rho _{i+1} - \rho _i}{\sigma (\delta )} \leq 2
\qquad \quad \hbox{for $i=1, \dots , m-1$.}
\end{equation}
Let
$$R_\delta ^i = \{(\rho , \vartheta ) \in R_{\frac {\sigma (\delta )}{\delta},\delta} : \rho _i \leq
\rho \leq \rho _{i+1}\} \qquad \quad \hbox{for $i=1, \dots ,
m-1$.}$$
Define
$$\hat{u}(\rho ) = \frac 1{2\pi} \int _0^{2\pi} u(\rho , \vartheta
)d\vartheta\,\qquad \hbox{for a.e. $\rho \in (\overline \rho -
\varepsilon , \hat \rho - \varepsilon )$.}$$
 We have that
 \begin{equation}\label{ch33}
 \| |u|^p \|_{L^A(R_{\frac {\sigma (\delta )}{\delta},\delta} )}
 \leq  2^{p-1}\| |u-\hat{u} |^p \|_{L^A(R_{\frac {\sigma (\delta )}{\delta},\delta})} +
2^{p-1} \| |\hat{u}|^p \|_{L^A(R_{\frac {\sigma (\delta
)}{\delta},\delta} )}\,,
 \end{equation}
 where $\hat u$ is regarded as a function defined on $R_{\frac {\sigma (\delta )}{\delta},\delta}$,
 and  $A$ is the Young function introduced in Step 2. Furthermore,
 \begin{equation}\label{ch34}
  \| |u-\hat{u} |^p \|_{L^A(R_{\frac {\sigma (\delta )}{\delta},\delta} )} \leq \sum _{i=1}^{m-1} \| |u-\hat{u} |^p \|_{L^A(R _\delta
  ^i)}\,.
 \end{equation}
The mean value of $u - \hat{ u}$ over each $R_\delta
^i$ is $0$.
The manifolds $R_\delta ^2 , \cdots , R_\delta ^{m-2}$ agree, up
to translations, with the same cylinder. The manifolds $R_\delta
^1$ and $ R_\delta ^{m-1}$  also coincide, up to isometries.
Moreover,
$$\hdue (R_\delta ^i ) \approx \sigma (\delta )^2\,,$$
owing to  assumption \eqref{approx10}.
An analogous scaling argument as in the proof of Step 2 tells us
\begin{equation}\label{ch35}
\| |u-\hat{u} |^p \|_{L^A(R _\delta ^i)} \leq C \frac{\sigma
(\delta )^{p-2}}{A^{-1}(C\sigma (\delta )^{-2})} \|\nabla
u\|_{L^p(R _\delta
  ^i)}^p\,,
  \end{equation}
  for some constant $C$.
  Next, since
  \begin{equation*}
  |\hat{u} (\rho )|^p \leq C \bigg(|\hat{u} (\hat \rho -
  \varepsilon)|^p + \bigg(\int _{\overline \rho -\varepsilon }^{\hat \rho -\varepsilon
  }|\hat{u} '(r )|dr \bigg)^p\bigg)\, \qquad \hbox{for a.e. $\rho \in
  (\overline \rho -\varepsilon , \hat \rho -\varepsilon )$,}
  \end{equation*}
  one has that
  \begin{align}\label{ch37}
\| & |\hat{u}|^p \|_{L^A(R_{\frac {\sigma (\delta
)}{\delta},\delta})} \\ \nonumber & \leq C _0|\hat{u} (\hat \rho -
  \varepsilon)|^p \| 1 \|_{L^A(R_{\frac {\sigma (\delta )}{\delta},\delta}
 )} + C _0\bigg\|\bigg(
 \int _{\overline \rho -\varepsilon }^{\hat \rho -\varepsilon
  }    \frac 1{2\pi} \int _0^{2\pi} |u_\rho (\rho , \vartheta
)|d\vartheta d\rho\bigg)^p\bigg\|_{L^A(R_{\frac {\sigma (\delta
)}{\delta},\delta}
 )}
 \\ \nonumber &
 \leq \frac{C_0}{A^{-1}(1/\hdue (R_{\frac {\sigma (\delta )}{\delta},\delta} ))} |\hat{u} (\hat
\rho -
  \varepsilon)|^p  + C_0 \bigg\|\bigg(
 \int _{\overline \rho -\varepsilon }^{\hat \rho -\varepsilon
  }    \frac 1{2\pi} \int _0^{2\pi} |\nabla u|d\vartheta d\rho\bigg)^p\bigg\|_{L^A(R_{\frac {\sigma (\delta )}{\delta},\delta}
 )}
 \\ \nonumber &
 \leq \frac{C_0}{A^{-1}(1/\hdue (R_{\frac {\sigma (\delta )}{\delta},\delta} ))} |\hat{u} (\hat
\rho -
  \varepsilon)|^p
\\ \nonumber
 & \qquad + \frac{C_0}{\min _{\rho \in (\overline \rho -\varepsilon , \hat \rho -\varepsilon )} \phi (\rho )^p} \bigg\|\bigg(
 \int _{\overline \rho -\varepsilon }^{\hat \rho -\varepsilon
  }    \frac 1{2\pi} \int _0^{2\pi} |\nabla u|\phi (\rho )d\vartheta d\rho\bigg)^p\bigg\|_{L^A(R_{\frac {\sigma (\delta )}{\delta},\delta}
 )}
\\ \nonumber
 &
 \leq \frac{C_0}{A^{-1}(1/\hdue (R_{\frac {\sigma (\delta )}{\delta},\delta} ))} |\hat{u} (\hat
\rho -
  \varepsilon)|^p  + \frac{C_1}{\sigma (\delta )^p} \bigg\|\bigg(
 \int _{R_{\frac {\sigma (\delta )}{\delta},\delta}} |\nabla u|d\hdue \bigg)^p\bigg\|_{L^A(R_{\frac {\sigma (\delta )}{\delta},\delta}
 )}
\\ \nonumber &
 \leq \frac{C_0}{A^{-1}(1/\hdue (R_{\frac {\sigma (\delta )}{\delta},\delta} ))} |\hat{u} (\hat
\rho -
  \varepsilon)|^p  + \frac{C_1}{\sigma (\delta )^p} \hdue (R_{\frac {\sigma (\delta )}{\delta},\delta} )^{p-1}
  \|\nabla u\|_{L^p(R_{\frac {\sigma (\delta )}{\delta},\delta}
 )}^p\|1\|_{L^A(R_{\frac {\sigma (\delta )}{\delta},\delta}
 )}
\\ \nonumber &
 \leq \frac{C_2}{A^{-1}(C_3/(\delta \sigma (\delta )))} \bigg(\frac
 1{\sigma (\delta )} \int _{\partial N_{\frac {\sigma (\delta )}{\delta},\delta}} |u| d\huno \bigg)^p
  + \frac{C_2\delta ^{p-1}\sigma (\delta )^{p-1}}{\sigma (\delta )^p  A^{-1}(C_3/\delta \sigma (\delta )) }
   \|\nabla u\|_{L^p(R_{\frac {\sigma (\delta )}{\delta},\delta}
 )}^p
\\ \nonumber &
 \leq \frac{C_4}{\sigma (\delta )A^{-1}(C_3/(\delta \sigma (\delta )))} \int _{\partial N_{\frac {\sigma (\delta )}{\delta},\delta}} |u|^p d\huno
  + \frac{C_2\delta ^{p-1}}{\sigma (\delta )  A^{-1}(C_3/(\delta \sigma (\delta ))) }  \|\nabla u\|_{L^p(R_{\frac {\sigma (\delta )}{\delta},\delta}
 )}^p
 \,,
 \end{align}
 for suitable constants $C_i$, $i=0, \dots ,4$.
Here, we have made use of the fact that
$$\hdue (R_{\frac {\sigma (\delta )}{\delta},\delta}  ) \approx \delta \sigma (\delta
)\,.$$ An approximation argument for $u$ by smooth functions is
also required.
Owing to \eqref{D}, for any $C>0$, there exists a constant $C'>0$ such that
\begin{equation}\label{ch38}
\frac{\sigma (\delta) ^{p-2}}{  A^{-1}(C/ \sigma (\delta
)^2) } \leq \frac{C'\delta ^{p-1}}{\sigma (\delta )
A^{-1}(C'/(\delta \sigma (\delta ))) }\,.
\end{equation}
Thus, from \eqref{ch33}--\eqref{ch38} one deduces that there exists a constant $C>0$ such that
\begin{multline}\label{ch39}
\| |u|^p \|_{L^A(R_{\frac {\sigma (\delta )}{\delta},\delta}
 )} \\ \leq \frac{C}{\sigma (\delta )A^{-1}(C/(\delta \sigma (\delta )))} \int _{\partial N_{\frac {\sigma (\delta )}{\delta},\delta}} |u|^p d\huno
  + \frac{C\delta ^{p-1}}{\sigma (\delta )  A^{-1}(C/(\delta \sigma (\delta ))) }
  \int _{R_{\frac {\sigma (\delta )}{\delta},\delta} }|\nabla u|^p d\hdue\,.
 \end{multline}
Consequently,
\begin{align}\label{ch40}
\| & |u|^p \|_{L^A(\cup _k R^k
 )}  \leq \sum _{k \in \N} \| |u|^p \|_{L^A( R ^k)}
 \\ \nonumber & \leq
\sum _{k \in \N} \frac{C}{\sigma (\delta _k)A^{-1}(C/(\delta
_k\sigma (\delta _k)))} \int _{\partial N^k} |u|^p d\huno
  + \sum _{k \in \N} \frac{C\delta _k^{p-1}}{\sigma (\delta _k)  A^{-1}(C/(\delta _k\sigma (\delta _k))) }
  \int _{R^k}|\nabla u|^p d\hdue\,.
 \end{align}
By
\eqref{D}  and \eqref{ch17},
\begin{equation}\label{ch41}
\frac{\delta _k^{p-1}}{\sigma (\delta _k)  A^{-1}(C/(\delta
_k\sigma (\delta _k))) } \leq \frac{\delta _k^{p-1}}{\sigma
(\delta _k)  A^{-1}(C'/\delta _k ^2)} \leq C''\,,
\end{equation}
for some positive constants $C'$ and $C''$. Thus,
\begin{align}\label{ch42}
\| |u|^p \|_{L^A(\cup _k R ^k)}
 \leq
 \sum _{k \in \N} \frac{C}{\delta _k ^{p-1}} \int _{\partial N^k} |u|^p d\huno + C
 \|\nabla u\|_{L^p(\cup _k R^k)}^p
\,,
\end{align}
for some constant $C$. Hence, since $\sigma (\delta _k)$ is
bounded for $k \in \Z$, we deduce from \eqref{ch31'} that
\begin{equation}\label{ch43}
\| |u|^p \|_{L^A(\cup _k R ^k)}^{1/p} \leq C \big(\|\nabla
u\|_{L^p(M)} + \| u\|_{L^p(M)}\big).
\end{equation}
\par\noindent
{\bf Step 5}.  A variant of \cite[Theorem 2.3.2]{Ma2}, with
analogous proof, tells us that given a (2-dimensional) Riemannian
manifold $Z$ with $\mathcal H ^2 (Z) < \infty$, and a Young
function $B$, the inequality
\begin{equation}\label{ch44}
\| |u|^p \|_{L^B(Z)}^{1/p} \leq C \big(\|\nabla u\|_{L^p(Z)} + \|
u\|_{L^p(Z)}\big)
\end{equation}
 holds for some constant $C$ and for every $u \in W^{1,p}(Z)$ if and only if
\begin{equation}\label{ch45}
\frac{1}{B^{-1}(1/s)} \leq C' \nu _{Z,p}(s) \qquad \hbox{for $s \in (0,
\hdue (Z)/2)$},
\end{equation}
 for some constant $C'$.
In Step 3 we have observed that a regular neighbourhood of $\cup
_k
\partial N^k$ is a planar domain fulfilling the cone
property. Hence, the standard Sobolev inequality holds on $M
\setminus (\cup _k N^k)$, and, consequently, \eqref{ch45} holds
with $Z=M \setminus (\cup _k N^k)$ and
 $B(t) = t^{\frac{2}{2-p}}$ if $1 \leq p <2$, and with $B(t)=t^a$
for any $a \geq 1$ if $p=2$. Thus, since the right-hand side of
\eqref{ch18} is equivalent to a non-decreasing function,
inequality \eqref{ch45} also holds with $B=A$. Hence, there exists
a constant $C$ such that
%
%
\begin{equation}\label{ch46}
\| |u|^p \|_{L^A(M \setminus (\cup _k N^k))}^{1/p} \leq C
\big(\|\nabla u\|_{L^p(M \setminus (\cup _k N^k))} + \| u\|_{L^p(M
\setminus (\cup _k N^k))}\big)\,
\end{equation}
for $u \in W^{1,p}(M)$. Combining \eqref{ch32}, \eqref{ch43} and
\eqref{ch46} tells us that
\begin{equation}\label{ch47}
\| |u|^p \|_{L^A(M )}^{1/p} \leq C \big(\|\nabla u\|_{L^p(M )} +
\| u\|_{L^p(M )}\big)
\end{equation}
for some constant $C$ and for every $u \in W^{1,p}(M)$.
Hence,
\begin{equation}\label{ch48}
\frac{1}{A^{-1}(1/s)} \leq C \nu _{M,p}(s) \qquad \hbox{for $s \in
(0, \hdue (M)/2)$}
\end{equation}
and  \eqref{ch4} follows, owing to \eqref{ch17}.
\par\noindent
{\bf Part II}. Here we show that, if $p \geq 1$ and $\sigma$ is
non-decreasing and is of class $ \Delta _2$ near $0$, then
inequality \eqref{ch-3} holds.
%
%
%
%
Consider the sequence of condensers $(Q^k\cup P^k, N^k)$. Let
$\{u_k\}$ be the sequence of Lipschitz continuous functions given
by $u_k = 1$ in $Q^k\cup P^k$, $u_k = 0$ in $M \setminus N^k$ and
such that $u_k$ depends only on $\rho$ and is a linear function of
$\rho$ in $R^k$. For $k \in \N$, we have that
\begin{equation}\label{ch50}
\hdue (Q^k\cup P^k) \approx \delta _k^2\,,
\end{equation}
and
\begin{equation}\label{ch51}
\int _M |\nabla u_k|^p d \hdue \approx \frac{\hdue (R^k)}{\delta
_k^p} \approx \frac {\sigma (\delta _k)}{\delta _k^{p-1}}\,.
\end{equation}
 Thus, there exist constants $C$ and $C'$ such that
\begin{equation}\label{ch53}
\nu _{M,p} (C\delta _k^2) \leq C_p(Q^k\cup P^k, N^k) \leq \frac
{C'\sigma (\delta _k)}{\delta _k^{p-1}}\,.
\end{equation}
It is easily seen that \eqref{ch53} continues to hold with $\delta
_k$ replaced by any $s \in (0, \m2 )$. Hence \eqref{ch-3} follows.
\par\noindent
The proof is complete. \qed

\medskip
\par\noindent
{\bf Acknowledgements}. This research was partially supported by
the research project of MIUR ``Partial differential equations and
functional inequalities: quantitative aspects, geometric and
qualitative properties, applications", by the Italian research
project ``Geometric properties of solutions to variational
problems" of GNAMPA (INdAM) 2006, and by the UK and Engineering
and Physical Sciences Research Council via the grant EP/F005563/1.



\begin{thebibliography}{99}



\bibitem[AFP]{AFP} L. Ambrosio, N. Fusco and D. Pallara,  Functions of bounded
variation and free discontinuity problems,  Oxford University
Press, Oxford, 2000.



\bibitem[BGM]{BGM} M. Berger, P. Gauduchon and E. Mazet,  ``Le spectre d'une vari\'et\'e Riemannienne", Lecture notes in Mathematics 194,
Springer-Verlag, Berlin, 1971.

\bibitem[BC]{BC} I. Benjamini and J. Cao, A new isoperimetric theorem for surfaces
of variable curvature, \emph{Duke Math. J.} {\bf 85} (1996),
359–-396.
%

\bibitem[Bou]{Bou} J. Bourgain, Geodesic restrictions and $L^p$-estimates for eigenfunctions of Riemannian surfaces,
 \emph{Amer. Math. Soc. Tranl.}  {\bf 226} (2009),
27--25.



\bibitem[Br]{Br} B. Brooks, The bottom of the spectrum of a Riemannian covering, \emph{J. Reine Angew. Math.} {\bf357} (1985),
101–-114.

\bibitem[BuZa]{BuZa} Yu. D. Burago and V. A. Zalgaller, ``Geometric inequalities", Springer-Verlag,
Berlin, 1988.

\bibitem[BD]{BD}
V. I. Burenkov and E. B. Davies,  Spectral stability of the
Neumann Laplacian, \emph{J. Diff. Eq.} {\bf 186} (2002), 485--508.

\bibitem[Cha]{Chavel}  I. Chavel, ``Eigenvalues in Riemannian geometry",  Academic Press,
New York,  1984.

\bibitem[CF]{CF} I. Chavel and E. A. Feldman, Modified isoperimetric
constants, and large time heat diffusion in Riemannian manifolds,
\emph{Duke Math. J.} { \bf 64} (1991), 473--499.

\bibitem[Che]{Cheeger} J. Cheeger, A lower bound for the smallest eigevalue of the Laplacian, in \emph{Problems in analysis},
195--199, Princeton Univ. Press, Princeton, 1970.



\bibitem[CGY]{CGY} F. Chung, A. Grigor'yan and S.-T. Yau, Higher eigenvalues and
isoperimetric inequalities on Riemannian manifolds and graphs,
\emph{Comm. Anal. Geom.} { \bf 8} (2000),  969--1026.
\bibitem[Ci1]{Crelative} A. Cianchi, On relative isoperimetric inequalities
in the plane, \emph{Boll. Un. Mat. Ital.}  {\bf 3-B} (1989),
289--326.
\bibitem[Ci2]{Cpoincare} A. Cianchi, O sharp form of Poincar\'e type inequalities on balls and spheres, \emph{Z. Angew. Math. Phys. (ZAMP)}
 {\bf 40} (1989), 558--569.

\bibitem[Ci3]{Celliptic} A. Cianchi, Elliptic equations on manifolds and
isoperimetric inequalities, \emph{Proc. Royal Soc. Edinburgh} {\bf
114A} (1990), 213--227.

\bibitem[Ci4]{Cmoser} A. Cianchi, Moser-Trudinger inequalities without boundary conditions and isoperimetric problems,
 \emph{Indiana Univ. Math. J.} {\bf
54} (2005), 669--705.
\bibitem[CEG]{CEG} A. Cianchi, D. E. Edmunds and  P.Gurka, On weighted
Poincar\'e inequalities, \emph{Math. Nachr.} {\bf 180} (1996),
 15--41.

\bibitem[CM]{CM} A. Cianchi and  V. G. Maz'ya, On the discreteness of the spectrum of the Laplacian on complete Riemannian manifolds,
 \emph{J. Diff. Geom.}, to appear.

\bibitem[CGL]{CGL} T. Coulhon, A. Grigor'yan and  D. Levin, On isoperimetric profiles of
product spaces, \emph{Comm. Anal. Geom.} {\bf 11} (2003), 85--120.

\bibitem[CH]{CH} R. Courant  and D. Hilbert, ``Methoden der mathematischen Physik", Springer, Berlin,
1937.


\bibitem[Da]{Da} G.Dal Maso, Notes on capacity theory, manuscript.



\bibitem[DS]{DS} E. B. Davies and  B. Simon, Spectral properties of the Neumann Laplacian of horns, \emph{Geom.
Funct. Anal.} {\bf 2} (1992), 105--117


\bibitem[Do1]{Do1} H. Donnelly, Bounds for eigefunctions of the Laplacian
on compact Riemannian manifolds, \emph{J. Funct. Anal.} {\bf 187}
(2001), 247--261.

\bibitem[Do2]{Do2} H. Donnelly, Eigenvalue estimates for certain noncompact manifolds, \emph{Michigan Math. J.}
{\bf 31} (1984), 349--357.


\bibitem[Es]{Escobar} J. F. Escobar, On the spectrum of the Laplacian on complete Riemannian manifolds, \emph{Comm. Part. Diff. Equat.} {\bf 11}
(1986), 63--85.

%


\bibitem[Ga]{Gallot} S. Gallot, In\'egalit\'es isop\'erim\'etriques et analitiques sur les vari\'et\'es riemanniennes, \emph{Asterisque}
{\bf 163} (1988), 31--91.


\bibitem[Gr1]{Gr1} A. Grigor'yan, On the existence of positive fundamental solution of the Laplace
equation on Riemannian manifolds,
 \emph{Mat.. Sbornik} {\bf 128} (1985), 354--363 (Russian); English translation: \emph{Math. USSR Sb.}  {\bf 56} (1987),
 349--358.

\bibitem[Gr2]{Gr} A. Grigor'yan,  Isoperimetric inequalities and capacities on Riemannian
manifolds, in The Maz'ya anniversary collection, Vol. 1 (Rostock,
1998), 139--153, \emph{Oper. Theory Adv. Appl.}, 109, Birkhäuser,
Basel, 1999.

\bibitem[GP]{GP}R. Grimaldi and  P. Pansu,  Calibrations and isoperimetric profiles,
\emph{ Amer. J. Math.} {\bf 129} (2007),  315--350.





\bibitem[HK]{HK} P. Hai\l asz and  P. Koskela, Isoperimetric inequalites and imbedding theorems in irregular
domains, \emph{J. London Math. Soc.} {\bf 58} (1998), 425--450.

\bibitem[HSS]{HSS} R. Hempel, L. Seco and  B. Simon,  The essential spectrum of Neumann Laplacians on some bounded singular domains,
\emph{J. Funct. Anal.} {\bf 102} (1991), 448--483.



\bibitem[Hi]{Hille} E. Hille, Non-oscillatory theorems, \emph{Trans. Amer. math. Soc.} {\bf 64} (1948), 234--252.


\bibitem[HHN]{HHN} M. Hoffmann-Ostenhof, T. Hoffmann-Ostenhof, N. Nadirashvili, On
the multiplicity of eigenvalues of the Laplacian on surfaces,
\emph{Ann. Global Anal. Geom.} {\bf 17} (1999),  43--48.

\bibitem[JMS]{JMS} V. Jaksic, S. Molchanov and  B. Simon,  Eigenvalue asymptotics of the Neumann
Laplacian of regions and manifolds with cusps, \emph{J. Funct.
Anal.} {\bf 106} (1992), 59--79.



\bibitem[KK]{KK} I. S. Kac and  M. G. Krein, Criterion for discreteness of the spectrum of a singular string,
\emph{Izv. Vyss. Uchebn. Zaved Mat.} {\bf 2} (1958), 136--153
(Russian).



\bibitem[KM]{KM} T. Kilpel\"ainen and  J. Mal\'y,  Sobolev inequalities on sets with
irregular boundaries,  \emph{Z. Anal. Anwendungen} {\bf 19}
(2000), 369--380.

\bibitem[Kl]{Kl} B. Kleiner, An isoperimetric comparison theorem, \emph{Invent. Math.} {\bf 108} (1992), 37--47.

\bibitem[La]{La} D. A. Labutin, Embedding of Sobolev spaces on H\"older domains,
\emph{Proc. Steklov Inst. Math.} \textbf{227} (1999), 163--172
(Russian); English translation: {\em Trudy Mat. Inst.} {\bf 227}
(1999), 170--179.

\bibitem[LP]{LP} P.-L. Lions and  F. Pacella, Isoperimetric inequalities for convex cones,
\emph{Proc. Amer. Math. Soc.} \textbf{109} (1990), 477--485.

\bibitem[MZ]{MZ} J. Mal\'y and  W. P. Ziemer, ``Fine
regularity of solutions of elliptic partial differential
equations",
American Mathematical Society, Providence, 1997.
\bibitem[Ma1]{Ma0} V. G. Maz'ya, Classes of regions and imbedding theorems for function spaces,
{\em Dokl. Akad. Nauk. SSSR} {\bf 133} (1960), 527--530 (Russian);
English translation: {\em Soviet Math. Dokl.} {\bf 1} (1960),
882--885.
\bibitem[Ma2]{Ma0bis} V. G. Maz'ya, Some estimates of solutions of second-order elliptic equations,
{\em Dokl. Akad. Nauk. SSSR} {\bf 137} (1961), 1057--1059
(Russian); English translation: {\em Soviet Math. Dokl.} {\bf 2}
(1961), 413--415.
\bibitem[Ma3]{Ma4} V. G. Maz'ya, On p-conductivity and theorems on embedding certain functional
spaces into a C-space,  \emph{Dokl. Akad. Nauk SSSR} 140 (1961),
299--302 (Russian).
\bibitem[Ma4]{Ma5} V. G. Maz'ya,  On the solvability of the Neumann problem,
\emph{ Dokl. Akad. Nauk SSSR} 147 (1962), 294--296 (Russian).
\bibitem[Ma5]{Ma6} V. G. Maz'ya,  The Neumann problem in regions with nonregular
boundaries,  \emph{Sibirsk. Mat. \v{Z}.} 9 (1968), 1322--1350
(Russian).
\bibitem[Ma6]{Ma3} V. G. Maz'ya, On weak solutions of the Dirichlet and Neumann problems,
{\em Trusdy Moskov. Mat. Ob\v s\v c.} {\bf 20} (1969), 137--172
(Russian); English translation: \emph{Trans. Moscow Math. Soc.}
{\bf 20} (1969), 135--172.
\bibitem[Ma7]{Ma2} V. G. Maz'ya,  ``Sobolev spaces", Springer-Verlag, Berlin, 1985.


\bibitem[MHH]{HHM} F. Morgan, H. Howards and   M. Hutchings, The isoperimetric problem on surfaces of revolution of decreasing Gauss curvature,
\emph{Trans. Amer. Math. Soc.} {\bf 352} (2000), 4889--4909.

\bibitem[MJ]{MJ} F. Morgan, and   D. L. Johnson, Some sharp isoperimetric theorems for Riemannian manifolds,
\emph{Indiana Univ. Math. J.} {\bf 49} (2000), 1017--1041.

\bibitem[Na]{Na} N. Nadirashvili, Isoperimetric inequality for the second
eigenvalue of a sphere, \emph{J. Diff. Geom.} {\bf 61} (2002),
335--340.

\bibitem[On]{Oneil} R. O'Neil, Fractional intergration in Orlicz spaces, \emph{Trans. Amer. Math. Soc.} {\bf 115}
(1965), 300--328.

\bibitem[Pi]{Pi} Ch. Pittet, The isoperimetric profile of homogeneous Riemannian manifolds, \emph{J. Diff. Geom.} {\bf  54}
 (2000), 255--302.

\bibitem[Ri]{Ri} M. Ritor\'e, Constant geodesic curvature curves and isoperimetric
domains in rotationally symmetric surfaces, \emph{Comm. Anal.
Geom.} {\bf 9} (2001), 1093--1138.

\bibitem[Sa]{Saloff} L. Saloff-Coste, Sobolev inequalities in familiar and unfamiliar settings, in \emph{Sobolev spaces in mathematics},
Vol I, Sobolev type inequalities,  V.G.Maz'ya editor, Springer,
2009.



\bibitem[SS]{SS}H. F. Smith and  C. D. Sogge,  On the $L\sp p$ norm of spectral clusters for
compact manifolds with boundary, \emph{Acta Math.} {\bf 198}
(2007), 107--153.
\bibitem[So]{So} C. D. Sogge, Lectures on eigenfunctions of the
Laplacian, Topics in mathematical analysis, 337--360, Ser. Anal.
Appl. Comput., 3, World Sci. Publ., Hackensack, NJ, 2008.
\bibitem[SZ]{SZ} C. D. Sogge and  S. Zelditch, Riemannian manifolds with maximal
eigenfunction growth, \emph{Duke Math. J.} {\bf 114} (2002),
387--437.
\bibitem[Ya]{Ya} S. T. Yau, Isoperimetric constants and the first eigenvalue of a compact
manifold, \emph{Ann. Sci. E´cole Norm. Sup.} {\bf 8} (1975),
487-–507.
\end{thebibliography}
\end{document}